\documentclass{article}

\usepackage{amsmath,amsfonts, latexsym,graphicx}

\newcommand{\U}{{\mathcal U}}
\newcommand{\W}{{\mathcal W}}
\newcommand{\0}{{\mathbf 0}}
\newcommand{\C}{{\mathbb C}}
\newcommand{\Z}{{\mathbb Z}}
\newcommand{\R}{{\mathbb R}}
\newcommand{\Q}{{\mathbb Q}}
\newcommand{\D}{{\mathbb D}}
\newcommand{\N}{{\mathbb N}}
\newcommand{\cL}{{\mathbb L}}
\newcommand{\Proj}{{\mathbb P}}
\newcommand{\hyp}{{\mathbb H}}
\newcommand{\supp}{\operatorname{supp}}

\newcommand{\im}{\mathop{\rm im}\nolimits}

\newcommand{\Adot}{\mathbf A^\bullet}
\newcommand{\Fdot}{\mathbf F^\bullet}
\newcommand{\Gdot}{\mathbf G^\bullet}

\newcommand{\strat}{{\mathfrak S}}

\newcommand{\gecc}{{\operatorname{gecc}}}
\newcommand{\blfootnote}{\xdef\@thefnmark{}\@footnotetext}

\newtheorem{defn0}{Definition}[section]
\newtheorem{prop0}[defn0]{Proposition}
\newtheorem{conj0}[defn0]{Conjecture}
\newtheorem{thm0}[defn0]{Theorem}
\newtheorem{lem0}[defn0]{Lemma}
\newtheorem{corollary0}[defn0]{Corollary}
\newtheorem{example0}[defn0]{Example}
\newtheorem{remark0}[defn0]{Remark}
\newtheorem{question0}[defn0]{Question}

\newenvironment{defn}{\begin{defn0}\hskip -.06in .}{\end{defn0}}
\newenvironment{prop}{\begin{prop0}\hskip -.06in .}{\end{prop0}}

\newenvironment{thm}{\begin{thm0}\hskip -.06in .}{\end{thm0}}
\newenvironment{lem}{\begin{lem0}\hskip -.06in .}{\end{lem0}}
\newenvironment{cor}{\begin{corollary0}\hskip -.06in .}{\end{corollary0}}
\newenvironment{exm}{\begin{example0}\hskip -.06in .\rm}{\end{example0}}
\newenvironment{rem}{\begin{remark0}\hskip -.06in .\rm}{\end{remark0}}

\newcommand{\defref}[1]{Definition~\ref{#1}}
\newcommand{\propref}[1]{Proposition~\ref{#1}}
\newcommand{\thmref}[1]{Theorem~\ref{#1}}
\newcommand{\lemref}[1]{Lemma~\ref{#1}}
\newcommand{\corref}[1]{Corollary~\ref{#1}}
\newcommand{\exref}[1]{Example~\ref{#1}}
\newcommand{\secref}[1]{Section~\ref{#1}}
\newcommand{\remref}[1]{Remark~\ref{#1}}

\newcommand{\qed}{\mbox{$\Box$}}
\newenvironment{proof}{\noindent {\bf Proof.}}{\qed\vskip 6pt}

\title{Enriched Relative Polar Curves and Discriminants}

\author{David B. Massey}

\date{}

\begin{document}

\baselineskip = 14pt

\maketitle

\bigskip

\hfill{\it In honor of L\^e D\~ung Tr\'ang on his 60th birthday}

\bigskip

\begin{abstract}Let $(f, g)$ be a pair of complex analytic functions on a singular analytic space $X$. We give ``the correct'' definition of the relative polar curve of $(f, g)$, and we give a very formal generalization of L\^e's attaching result, which relates the relative polar curve to the relative cohomology of the Milnor fiber modulo a hyperplane slice. We also give the technical arguments which allow one to work with a derived category version of the discriminant and Cerf diagram of a pair of functions. From this, we derive a number of generalizations of results which are classically proved using the discriminant. In particular, we give applications to families of isolated ``critical points''.
\end{abstract}

\sloppy

%\newpage

%\tableofcontents

%\newpage

\section{L\^e's Attaching Result and Our Previous Generalization}\label{sec:genle}

\footnotetext[1]{The author would like to thank Tsukuba University, Tokyo Science University, Hokkaido University, and K. Takeuchi, M. Oka, and T. Ohmoto; Section 6 of this paper was written during a visit to these universities with the support of these mathematicians.
\newline AMS subject classifications 32B15, 32C35, 32C18, 32B10.
\newline   keywords: polar curve, discriminant, Milnor fiber, nearby cycles, $a_f$ condition.}

Let  $\U$ be an open neighborhood of the origin in $\C^{n+1}$, and let ${\tilde f}:\U\rightarrow \C$ be a complex analytic function. We assume that $\mathbf 0\in V(\tilde f):=\tilde f^{-1}(0)$. We let $\Sigma \tilde f$ denote the critical locus of $\tilde f$.

\bigskip

In this paper, we describe an improvement/generalization of what is now a classic result in the study of singularities: the attaching result of L\^e in \cite{leattach}, which tells one how many $n$-cells are attached, up to homotopy, to a hyperplane slice of the Milnor fiber of $\tilde f$ in order to obtain the Milnor fiber, $F_{\tilde f, \mathbf 0}$, of $\tilde f$ itself.

\smallskip

However, first, we must discuss the relative polar curve.

\smallskip

Fix a point $p\in \U$. Let $z_0$ denote a generic linear form on $\C^{n+1}$, which, in fact, we take as the first coordinate function, after possibly performing a generic linear change of coordinates.

In  \cite{hammlezariski}, \cite{teissiercargese}, \cite{leattach}, \cite{letopuse}, Hamm, Teissier, and L\^e define and use the relative polar curve (of $\tilde f$ with respect to $z_0$), $\Gamma^1_{\tilde f, z_0}$, to prove a number of topological results related to the Milnor fiber of hypersurface singularities. We shall recall some definitions and results here. We should mention that there are a number of different characterizations of the relative polar, all of which agree when $z_0$ is sufficiently generic; below, we have selected what we consider the easiest way of describing the relative polar curve as a set, a scheme, and a cycle.

\smallskip

As a set, $\Gamma^1_{\tilde f, z_0}$ is the closure of the critical locus of $(\tilde f, z_0)$ minus the critical locus of $\tilde f$, i.e., $\Gamma^1_{\tilde f, z_0}$ equals $\overline{\Sigma(\tilde f, z_0)-\Sigma \tilde f}$, as a set. If $z_0$ is sufficiently generic for $\tilde f$ at $p$, then, in a neighborhood of $p$, $\Gamma^1_{\tilde f, z_0}$ will be purely one-dimensional (which includes the possibility of being empty); see \thmref{thm:genpolar} below. 

It is not difficult to give $\Gamma^1_{\tilde f, z_0}$ a scheme structure. We use $(z_0, \dots, z_n)$ as coordinates on $\U$. If $\Gamma^1_{\tilde f, z_0}$ is purely one-dimensional at $p$, then, at points $x$ near, but unequal to, $p$,  $\Gamma^1_{\tilde f, z_0}$ is given the structure of the scheme $\displaystyle V\left(\frac{\partial\tilde f}{\partial z_1}, \dots, \frac{\partial\tilde f}{\partial z_n}\right)$. One can also remove ``algebraically'' any embedded components of $\Gamma^1_{\tilde f, z_0}$ at $p$ by using {\it gap sheaves\/}; see Chapter 1 of \cite{lecycles}.

In practice, all topological applications of the relative polar curve use only its structure as an analytic cycle (germ), that is, as a locally finite sum of irreducible analytic sets (or germs of sets) counted with integral multiplicities (which will all be non-negative). We remark here that these are cycles, {\bf not} cycle classes; we do not mean up to rational equivalence. The intersection theory that one needs here is the simple case of proper intersections inside smooth manifolds; see 8.2 of \cite{fulton} or our summary in Appendix A of \cite{numcontrol}. If $C$ is a one-dimensional irreducible germ of $\Gamma^1_{\tilde f, z_0}$ at $p$, and $x\in C$ is close to, but unequal to, $p$, then  the component $C$ appears in the cycle $\Gamma^1_{\tilde f, z_0}$ with multiplicity given by the Milnor number of $\tilde f_{|_H}$ at $x$, where $H$ is a generic affine hyperplane passing through $x$.

\smallskip The following theorem tells one the relative polar curve has nice properties for a generic choice of the linear form $z_0$.

\begin{thm}{\rm (Hamm-L\^e)}\label{thm:genpolar} For a generic choice of $z_0$, 

\begin{enumerate}
\item $\Gamma^1_{\tilde f, z_0}$ is purely one-dimensional at $p$;

\item $\Gamma^1_{\tilde f, z_0}$ properly intersects $V(\tilde f-\tilde f(p))$ at $p$, i.e., $p$ is an isolated point in $\Gamma^1_{\tilde f, z_0}\cap V(\tilde f-\tilde f(p))$;

\item the cycle $\Gamma^1_{\tilde f, z_0}$ is reduced (near $p$), i.e., each component through $p$ appears with multiplicity $1$.
\end{enumerate}
\end{thm}
\begin{proof} Items 1 and 2 are proved in 2.1 of \cite{hammlezariski}.  Item 3 is Lemma 2.2.1 of \cite{hammlezariski}. \end{proof}

\bigskip

Now, we can state the main result of \cite{leattach}.

\begin{thm} {\rm (L\^e, \cite{leattach})}\label{thm:leattach}  For generic $z_0$, up to homotopy, $F_{\tilde f, \mathbf 0}$ is obtained from $F_{{\tilde f}_{|_{V(z_0)}}, \mathbf 0}$ by attaching $\tau:=\big(\Gamma^1_{\tilde f, z_0}\cdot V(f)\big)_{\mathbf 0}$ $n$-cells.

In particular, $H^k(F_{\tilde f, \mathbf 0}, F_{{\tilde f}_{|_{V(z_0)}}, \mathbf 0}) =0$ if $k\neq n$, and $H^n(F_{\tilde f, \mathbf 0}, F_{{\tilde f}_{|_{V(z_0)}}, \mathbf 0}) \cong\mathbb Z^{\tau}.$
\end{thm}

\medskip

\begin{rem}\label{rem:cerf} It will be important for us to understand some of L\^e's set-up in \cite{leattach}. 

Let $B_\epsilon$ denote a closed $2n$-ball of radius $\epsilon$ centered at the origin in $\mathbb C^n$, and let $\mathbb D_\delta$ denote a closed disk of radius $\delta$ centered at the origin in $\mathbb C$.

Then, L\^e shows that, for $0<|\xi|\ll\delta\ll\epsilon\ll 1$, $(\mathbb D_\delta\times B_\epsilon)\cap \tilde f^{-1}(\xi)$ has the homotopy-type of the Milnor fiber of $\tilde f$ at $\mathbf 0$ and, of course,  $(\{0\}\times B_\epsilon)\cap  \tilde f^{-1}(\xi)$ is homeomorphic to the Milnor fiber, $F_{\tilde f_0, \mathbf 0}$, of $\tilde f_0:=\tilde f_{|_{V(z_0)}}$ at $\mathbf 0$. \thmref{thm:leattach} is obtained by applying Morse Theory to the map $|z_0|^2$ on $(\mathbb D_\delta\times B_\epsilon)\cap f^{-1}(\xi)$. 

 In \cite{leattach}, L\^e gives an extensive discussion of the discriminant and Cerf diagram of the map $G:=(z_0, \tilde f)$. The discriminant is $G(\Sigma G)$ and the Cerf diagram is $G(\Gamma^1_{\tilde f, z_0})$. One ``sees'' L\^e's attaching result graphically in the discriminant/Cerf diagram below. Down in the image of the map $(z_0, \tilde f)$, the pair $(F_{\tilde f, \mathbf 0}, F_{{\tilde f_0}, \mathbf 0})$ is represented by the pair $(L, \{a\})$, and one sees that the relative cohomology of $(F_{\tilde f, \mathbf 0}, F_{{\tilde f_0}, \mathbf 0})$ decomposes as a direct sum of local Morse data above each of the points where $L$ intersects the Cerf diagram, $C$, (the image of the relative polar curve). Now, the number of points in the intersection of $C$ and $L$, counted with multiplicities, is precisely $\tau=\big(\Gamma^1_{\tilde f, z_0}\cdot V(f)\big)_{\mathbf 0}$, which yields the theorem.

\begin{center}{\includegraphics[height=40mm, width =80mm]{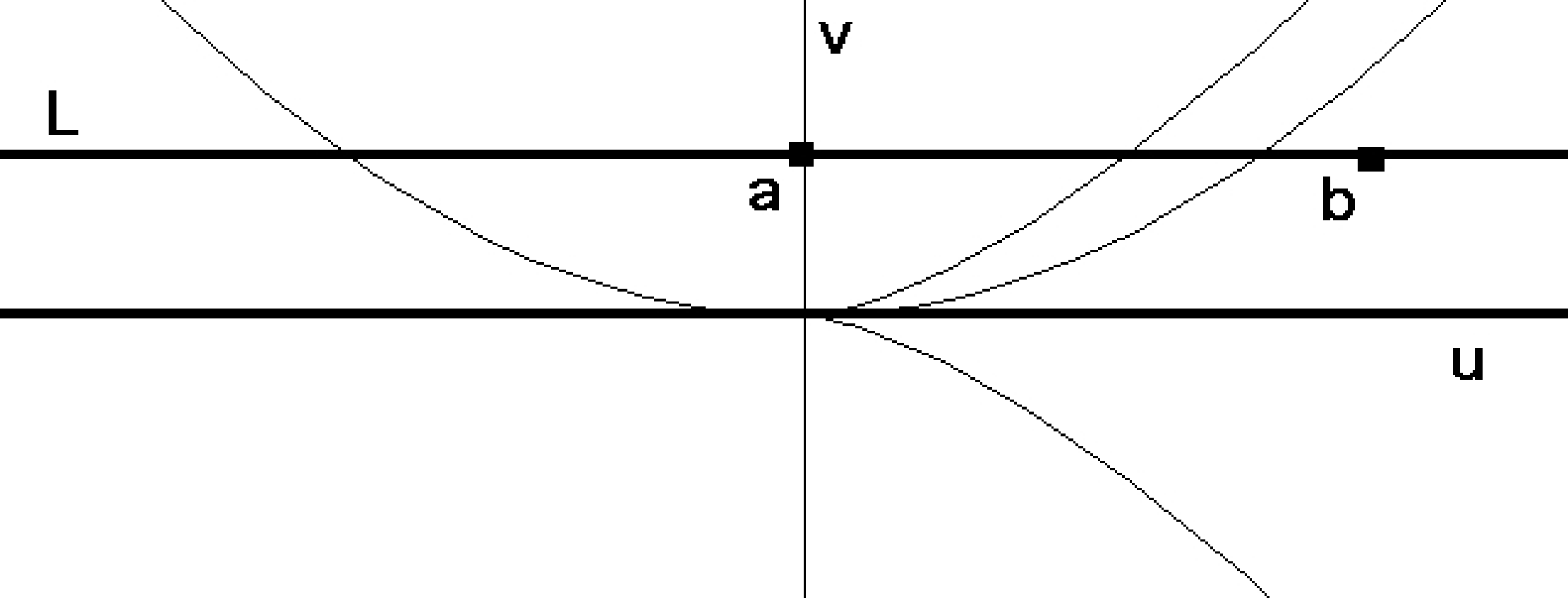}}\end{center}

From L\^e's discussion, it is clear that $F_{\tilde f_0, \mathbf 0}$ is also homeomorphic to $\widehat F_{\tilde f_0, \mathbf 0}:=(\{\nu\}\times B_\epsilon)\cap  \tilde f^{-1}(\xi)$, provided that $0<|\xi|\ll\nu<\delta\ll\epsilon\ll 1$; in fact, there is a homeomorphism from the pair $(F_{\tilde f, \mathbf 0}, F_{{\tilde f_0}, \mathbf 0})$ to $(F_{\tilde f, \mathbf 0}, \widehat F_{{\tilde f_0}, \mathbf 0})$  which induces an isomorphism between $H^*(F_{\tilde f, \mathbf 0}, F_{{\tilde f_0}, \mathbf 0})$ and $H^*(F_{\tilde f, \mathbf 0}, \widehat F_{{\tilde f_0}, \mathbf 0})$ and which induces the identity map on  $H^*(F_{\tilde f, \mathbf 0})$. In the Cerf diagram, the pair $(F_{\tilde f, \mathbf 0}, \widehat F_{{\tilde f_0}, \mathbf 0})$ is represented by $(L, \{b\})$ and, using an argument which is essentially the same as in the paragraph above, one concludes that $H^*(F_{\tilde f, \mathbf 0}, \widehat F_{{\tilde f_0}, \mathbf 0})\cong\Z^\tau$.

This is an important observation, because in the formalism of the derived category and vanishing cycles, $H^{k+1}(F_{\tilde f, \mathbf 0}, \widehat F_{{\tilde f_0}, \mathbf 0})$ is isomorphic to the stalk cohomology at the origin of the vanishing cycles along $z_0$ of the nearby cycles along $f$ of the constant sheaf on $\mathcal U$, i.e., $H^{k+1}(F_{\tilde f, \mathbf 0}, \widehat F_{{\tilde f_0}, \mathbf 0})\cong H^k(\phi_{z_0}\psi_f\mathbb Z^\bullet_{\mathcal U})_{\mathbf 0}$ (here, we do not distinguish between $z_0$ and ${z_0}_{|_{V(f)}}$). If we include the correct shifts, then we know that $\mathbb Z^\bullet_{\mathcal U}[n+1]$ is perverse, and that $\phi_{z_0}[-1]$ and $\psi_f[-1]$ take perverse sheaves to perverse sheaves; hence, we prefer to write $H^k(\phi_{z_0}[-1]\psi_f[-1]\mathbb Z^\bullet_{\mathcal U}[n+1])_{\mathbf 0}\cong H^{k+n}(F_{\tilde f, \mathbf 0}, \widehat F_{{\tilde f_0}, \mathbf 0})$.

Thus, the results of L\^e in \cite{leattach} tell one that $H^k(\phi_{z_0}[-1]\psi_f[-1]\mathbb Z^\bullet_{\mathcal U}[n+1])_{\mathbf 0}$ is zero, unless $k=0$, and $H^0(\phi_{z_0}[-1]\psi_f[-1]\mathbb Z^\bullet_{\mathcal U}[n+1])_{\mathbf 0}\cong \mathbb Z^{\tau}$. 
\end{rem}

\medskip

The above relation between \thmref{thm:leattach} and iterated vanishing and nearby cycles appeared explicitly  in the work of Sabbah in \cite{sabbahprox} and in our own work in  \cite{hypercohom}. 

\bigskip

\noindent{\bf Summary of our Old Results from \cite{hypercohom}}

\medskip

We now wish to describe one of our primary results from \cite{hypercohom}, which is a substantial generalization of \thmref{thm:leattach}, but first we need recall some of our previous definitions. One of our main goals in the current paper is to replace these old definitions with more natural ones.

\medskip

Let $X$ be a closed analytic subspace of $\mathcal U$, and let $f:=\tilde f_{|_{X}}$. 

\smallskip

A {\it good stratification of $X$ relative to $f$\/} is a complex analytic stratification $\strat$ of $X$ such that all of the strata $S\in\strat$ are connected, $V(f)$ is a union of strata,  ${\strat}^{o}:=\{S\in \strat \ | \ S\not\subseteq V(f)\}$ is a Whitney stratification of $X-V(f)$ and such that, for every pair of strata $(S_\alpha, S_\beta)$ such that  $S_\alpha\not\subseteq V(f)$
and $S_\beta\subseteq V(f)$, Thom's $a_f$ condition is satisfied. In our setting this is
equivalent to: if $\mathbf p\in S_\beta\subseteq V(f)$ and $\mathbf p_i\in
S_\alpha\not\subseteq V(f)$ are such that $\mathbf p_i\rightarrow \mathbf p$  and
$T_{\mathbf p_i}V\left(f_{|_{S_\alpha}} -  f_{|_{S_\alpha}}(\mathbf p_i)\right)$
converges to some $\mathcal T$, then $T_\mathbf p S_\beta \subseteq \mathcal T$; note that
this  implies that the pair $(S_\alpha, S_\beta)$ must satisfy Whitney's condition (a).
In a good  stratification, we call the strata which comprise $V(f)$ the {\it good
strata\/}; we refer to the other strata as {\it outside strata}. Note that we do {\bf not} require that Whitney's condition (b) hold along good strata, or even that Whitney's condition (a) holds between pairs of good strata.

\smallskip

Fix a good stratification $\strat$ for $X$ relative to $f$.

\smallskip

Let $\tilde g:(\mathcal U, \mathbf 0) \rightarrow (\mathbb C, 0)$ be another analytic function, and let $g:=\tilde g_{|_{X}}$. If $Y$ is an analytic
subset of $X$, we  define $\Gamma_{f,g}(Y)$ to be the closure in $X$ of the critical
locus of $\Phi_{|_{Y -  \Sigma Y - V(f)}}$.  This is called the {\it relative polar
variety of $Y$ with respect to $f$ and  $g$}

For each $S\in\strat$,
$\Gamma_{{}_{f,  g}}(S)$  is thus the closure of the critical locus of 
$(f, g)_{|_{S - V(f)}}$.  The union  $\cup_{S\in\strat}\Gamma_{{}_{f,
g}}(S)$
 is called {\it the relative polar variety of $f$ and $g$ with respect to $\strat$},
and we denote it by $\Gamma_{{}_{f, g}}(\strat)$ (or, simply, $\Gamma_{{}_{f, g}}$ if
the stratification is clear).  Note that if $S \subseteq
V(f)$, then $\Gamma_{f,g}(S) = \emptyset$.

For each stratum $S\in\strat$, we define the {\it symmetric relative polar variety of $S$ with respect to $f$ and 
$g$},  $\widetilde\Gamma_{f, g}(S)$, to be the closure in $X$ of the critical locus of $(f, g)_{|_{S - V(f) - V(g)}}$.  We also  define {\it the symmetric relative polar variety of $f$ and $g$ with respect to $\strat$},
$\widetilde\Gamma_{{}_{f, g}}(\strat)$, to be the union  $\bigcup_{S\in\strat}\
\widetilde\Gamma_{{}_{f, g}}(S)$.  We use the term ``symmetric'' since we
obviously have $\widetilde\Gamma_{f, g}(Y) = \widetilde\Gamma_{g,f}(Y)$.

\vskip .1in

In the special case where $f$ and $g$ are such that $\widetilde\Gamma_{{}_{f, g}}(\strat)$ is one-dimensional, we naturally
refer to the  symmetric polar  variety as the {\it  symmetric polar
curve\/} and emphasize the fact that it is one- dimensional by writing
$\widetilde\Gamma^1_{{}_{f,  g}}(\strat)$.  In this case, we wish to give the symmetric polar curve the structure of  a cycle (actually, a cycle germ at the origin), so we must attach some multiplicity to each  component of this curve.
 
 To do this, for each (one-dimensional) component, $\nu$, of $\widetilde\Gamma^1_{{}_{f,  g}}(\strat)$, let $S_\nu$
denote the stratum which contains $\nu - \mathbf 0$ near the origin.  If $S_\nu$ is
itself one-dimensional, we assign the multiplicity $1$ to $\nu$ (that  is, we
consider $\nu$ with its reduced structure). Now, to each $\nu$ for which  $S_\nu$ is
not one-dimensional, we assign the multiplicity given by the Milnor number of the 
map $g$ restricted to $S_\nu \cap V(f - f(\mathbf p))$ at any point $\mathbf p \in \nu -
\mathbf 0$ sufficiently close to the origin.  We use here that $S_\nu \cap V(f - f(\mathbf p))$
is a  manifold at $\mathbf p$, and that $g$ restricted to this set has an
isolated critical  point at $\mathbf p$ since $\nu$ is one-dimensional.

\smallskip

The function $g$ is {\it tractable at the origin with respect to a good stratification $\strat$ of $X$ relative to $f:(X,  \mathbf 0)\rightarrow (\mathbb C, 0)$\/} if and only if ${\operatorname{dim}}_\mathbf 0\widetilde\Gamma^1_{f,g}(\strat)\leq 1$ and, for all good strata  $S_\alpha$, $g_{|_{S_\alpha}}$ has no critical points in a neighborhood of the origin except, perhaps, at the origin itself. 

\smallskip

We say that $g$ is {\it decent\/}
with respect to $\strat$ relative to $f$  provided that $d_\mathbf p\tilde g$ is not a degenerate covector (see \cite{stratmorse}) at any stratified critical point of $g$ restricted to $F_{f, \mathbf 0}-V(g)$  (with respect to the induced
stratification on $F_{f, \mathbf 0}$). See Section 1 of \cite{hypercohom}. Note that this condition is automatic if $\strat$ has only one  stratum not contained in $V(f)\cup V(g)$.

\medskip

By combining Proposition 1.12 and 1.14 of \cite{hypercohom}, we obtain:

\smallskip

\begin{prop}
Let $\strat$ be a good stratification of $X$
for $f$ at the origin. Then, for a generic choice of linear forms, $l$, $l$ is decent and tractable with respect to $\strat$ relative to $f$.
\end{prop}

\medskip

For $S\in{\strat}^{o}$, we let $\N_S$ and $\cL_S$ denote, respectively, the normal slice and link of the stratum $S$; see \cite{stratmorse}. In  \cite{hypercohom}, we proved:

\medskip

\begin{thm}\label{thm:old} {\rm(\cite{hypercohom}, Theorem 4.2)} Suppose that $g$ is tractable relative to $f$ with
respect to a good stratification $\strat$ of $X$ at $\mathbf 0$.  Let
$d_S$ denote the dimension of $S\in\strat$.  Let $\Fdot$ be a bounded complex of sheaves of $\mathbb Z$-modules on $X-V(f)$,  constructible with respect to $\strat^o$.

\vskip .1in

Then, for all $i$, $\mathbb H^i(F_{f, \mathbf 0}, F_{f_{|_{V(g)}}, \mathbf 0}  ; \Fdot)$
is a direct summand of $H^{i-1}(\phi_g\psi_f\Fdot)_\mathbf 0$, and  there exist
integers $j_S$ such that  
$$\mathbb H^i(F_{f, \mathbf 0}, F_{f_{|_{V(g)}}, \mathbf 0}
 ; \Fdot) \ \cong \  \bigoplus_{S\in{\strat}^{o}} \big(\mathbb H^{i-
d_S +1}(\N_S, \cL_S  ; \Fdot)\big)^{j_S} , $$ where
$j_S \geq\big(\widetilde\Gamma^1_{{}_{f,g}}(S)\cdot V(f)\big)_\mathbf
0$, with equality if $g$ is decent relative to $f$.

\vskip .1in

Furthermore, if $\Gamma_{{}_{f,g}}(\strat)$ has no components contained in $V(g)$
(i.e., if $\Gamma_{{}_{f,g}}(\strat) =  \widetilde\Gamma_{{}_{f,g}}(\strat)$), then
$$ H^{i-1}(\phi_g\psi_f\Fdot)_\mathbf 0 \ \cong  \ \mathbb H^i(F_{f, \mathbf 0},
F_{f_{|_{V(g)}}, \mathbf 0}  ; \Fdot). $$
\end{thm}

\bigskip

\noindent{\bf Summary of the Results of this Paper}

\medskip

What are the problems with \thmref{thm:old}? There are several. One is that the hypotheses are difficult to check. Another, related, problem is that it is unclear to what extent the hypotheses are necessary for the conclusion. A third issue is that the definition ``relative symmetric polar curve'' seems rather ad hoc.

In this paper, we ``fix'' these problems. Let $\strat^o(\Fdot)$ be the set of strata of $\strat$ such that $\mathbb H^*(\N_S, \cL_S  ; \Fdot)\neq 0$ and $f_{|_S}$ is not constant.  For each $S\in\strat^o(\Fdot)$, we will define an (ordinary) cycle $\Gamma_{f, \tilde g}(S)$. Using these cycles,  we will define (\defref{def:polarcurve}) the {\it graded, enriched relative polar cycle, $\big(\Gamma_{f, \tilde g}(\Fdot)\big)^\bullet$}. In each degree $k\in\Z$, $\big(\Gamma_{f, \tilde g}(\Fdot)\big)^k$ is a formal, locally finite, sum of irreducible analytic subsets of $X$ multiplied by modules over a fixed base ring; see Section 2 of \cite{singenrich} and \secref{sec:main}. When the underlying set, $\big|\big(\Gamma_{f, \tilde g}(\Fdot)\big)^\bullet\big|=\bigcup_{S\in\strat^o(\Fdot)}\left|\Gamma_{f, \tilde g}(S)\right|$ is purely one-dimensional at a point $p\in X$, we say that {\it the relative polar curve of $f$, with respect to $g$, with coefficients in $\Fdot$, is defined at $p$}. A principal theme of this paper is that this definition of the relative polar curve is {\bf the} correct definition in results on the cohomology level.

The intersection product that we use throughout our work is a mild extension of the intersection theory, mentioned above, of properly intersecting cycles in a complex manifold; see Section 2 of \cite{singenrich} for the fundamental properties. We use $\odot$ to denote this enriched intersection product.

\bigskip

Using our results in \cite{singenrich}, and continuing with the notation from above, we will quickly prove our first main theorem. Below, and throughout this paper, we adopt the convention that the empty set has dimension $-\infty$, so that for analytic space $Z$ and a point $p$, the condition that $\dim_pZ\leq 0$ means that either $\dim_pZ=0$ or $p\not\in Z$.

\medskip

\noindent{\bf Main Theorem 1}. {\rm (\thmref{thm:main1})}  {\it In a neighborhood of the origin, 
$$\supp \phi_g[-1]\psi_f[-1]\Fdot = V(f)\cap \big|\Gamma_{f, \tilde g}(\Fdot)\big|\subseteq V(g),$$
and, when $\dim_\0 V(f)\cap \big|\Gamma_{f, \tilde g}(\Fdot)\big|\leq 0$,
$$
H^k(\phi_g[-1]\psi_f[-1]\Fdot)_\0\cong \big(\big(\Gamma_{f, \tilde g}(\Fdot)\big)^k\odot V(f)\big)_\0,
$$
i.e., 
$$H^{k-1}(B_\epsilon\cap f^{-1}(\xi), B_\epsilon\cap f^{-1}(\xi)\cap g^{-1}(\nu);\ \Fdot)\ \cong\ \bigoplus_{S\in \strat^o(\Fdot)} \big(\mathbb H^{k-d_S }(\N_S, \cL_S  ; \Fdot)\big)^{j_S} , $$ 
where $0<|\xi|\ll |\nu|<\ll\epsilon\ll 1$, and 
$j_S = \big(\Gamma^1_{{}_{f,g}}(S)\cdot V(f)\big)_\mathbf 0$.}

\vskip .3in

Our proof of the above theorem is elegant, and very short, given existing results. However, it is not as intuitive as the discriminant/Cerf diagram argument, nor does it allow us to prove a number of related results, as we did in Section 4 of \cite{hypercohom}.

Hence, in \secref{sec:discrim} of this paper, we will prove the necessary technical results to push-down the complex $\Fdot$, restricted to a suitable neighborhood, via the map $(g, f)$. This will give us a derived category version of the discriminant and Cerf diagram, in which the standard intuitive proofs work without modification; see \thmref{thm:main2}. Somewhat surprisingly, the hypothesis that we need is precisely that of Main Theorem 1.

We prove:

\noindent{\bf Corollary to Main Theorem 2}. {\rm (\corref{cor:main2})}  {\it Suppose that $\dim_\0 V(f)\cap \big|\Gamma_{f, \tilde g}(\Fdot)\big|\leq 0$. Let $\big(\widehat\Gamma_{f, \tilde g}(\Fdot)\big)^\bullet$ denote the components of $\big(\Gamma_{f, \tilde g}(\Fdot)\big)^\bullet$ which are not contained in $V(g)$. Then,
\begin{enumerate}
\item
$$
\hyp^{k-1}(F_{f, \mathbf 0}, F_{f_{|_{V(g)}}, \mathbf 0} ; \Fdot) \cong \big(\big(\widehat\Gamma_{f, \tilde g}(\Fdot)\big)^k\odot V(f)\big)_\0;
$$

\item
$$
\hyp^{k-1}(F_{g, \mathbf 0}, F_{g_{|_{V(f)}}, \mathbf 0} ; \Fdot) \cong \big(\big(\widehat\Gamma_{f, \tilde g}(\Fdot)\big)^k\odot V(g)\big)_\0; and
$$

\item
$$
H^k(\phi_f[-1]\psi_g[-1]\Fdot)_\0\ \cong\ \big(\big(\widehat\Gamma_{f, \tilde g}(\Fdot)\big)^k\odot V(g)\big)_\0\ \oplus\ H^k(\psi_g[-1]\phi_f[-1]\Fdot)_\0.
$$
\end{enumerate}
}

\bigskip

In \secref{sec:afappl}, we will combine the results of \secref{sec:discrim} with Corollary 3.9 of \cite{vanaf} in order to obtain a relation between Thom's $a_f$ condition and the graded, enriched polar curve. In \secref{sec:families}, we show how the main theorems allow us to prove a number of familiar-looking results on families with isolated critical points.

\medskip

We belatedly thank Marc Levine for a number of helpful discussions involving our enriched intersection theory. We also thank L\^e D\~ung Tr\'ang for some very helpful proofreading, and for a number of suggestions which improved the presentation.

\section{Basics of Enriched Cycles}\label{sec:enrcycle}

In this section, we will recall the basic definitions that one needs for using enriched cycles ; these definitions are taken from Section 2 of \cite{singenrich}. There are a number of results from \cite{singenrich} which will be used in the proof of the main theorem in \secref{sec:main}. While we will not restate the needed results from \cite{singenrich} in this paper, the background material in this section will enable the reader to make sense of the definition of the graded, enriched relative polar curve and the proof of the main theorem in \secref{sec:main}.

\vskip .2in

\begin{defn} An {\it enriched cycle}, $E$, in $X$ is a formal, locally finite sum $\sum_V E_V[V]$, where the
$V$'s are irreducible analytic subsets of $X$ and the
$E_V$'s are finitely-generated $R$-modules. We refer to the $V$'s as the {\it components\/} of $E$, and to $E_V$ as the {\it
$V$-component module of
$E$}. Two enriched cycles are considered the same provided that all of the component modules are isomorphic. The underlying set
of $E$ is
$|E|:=
\cup_{{}_{E_V\neq 0}}V$.

If $C = \sum n_V[V]$ is an ordinary positive cycle in $X$, i.e., all of the
$n_v$ are non-negative integers, then there is a corresponding enriched cycle $[C]^{\operatorname{enr}}$ in which the
$V$-component module is the free
$R$-module of rank $n_V$. If $R$ is an integral domain, so that rank of an $R$-module is well-defined, then an enriched cycle
$E$ yields an ordinary cycle $[E]^{\operatorname{ord}}:=\sum_V (\operatorname{rk}(E_V))[V]$.

If $q$ is a finitely-generated module and $E$ is an enriched cycle, then we let
$qE:=\sum_V(q\otimes E_V)[V]$; thus, if $R$ is an integral domain and $E$ is an enriched cycle,
$[qE]^{\operatorname{ord}}=(\operatorname{rk}(q))[E]^{\operatorname{ord}}$ and if
$C$ is an ordinary positive cycle and $n$ is a positive integer, then $[nC]^{\operatorname{enr}}=R^n[C]^{\operatorname{enr}}$.

\vskip .1in

The (direct) sum of two enriched cycles $D$ and $E$ is given by $(D + E)_V := D_V\oplus E_V$. 

\vskip .1in

There is a partial ordering on enriched cycles given by: $D\leq E$ if and only if there exists an enriched cycle $P$ such that $D+P=E$. This relation is clearly reflexive and transitive; moreover, anti-symmetry follows from the fact that if $M$ and $N$ are Noetherian modules such that $M\oplus N\cong M$, then $N=0$.

\vskip .1in

If two irreducible analytic subsets $V$ and
$W$ intersect properly in $\U$, then the (ordinary) intersection cycle $[V]\cdot[W]$ is a well-defined positive cycle; we
define the enriched intersection product of $[V]^{\operatorname{enr}}$ and $[W]^{\operatorname{enr}}$ by
$[V]^{\operatorname{enr}}\odot[W]^{\operatorname{enr}} = ([V]\cdot[W])^{\operatorname{enr}}$. If $D$ and $E$ are enriched
cycles, and every component of $D$ properly intersects every component of $E$ in $\U$, then we say that {\it $D$ and $E$
intersect properly\/} in $\U$ and we extend the intersection product linearly, i.e., if $D=\sum_V D_V[V]$ and $E=\sum_W
E_W[W]$, then
$$ D\odot E:= \sum_{V, W} (D_V\otimes E_W)([V]\cdot [W])^ {\operatorname{enr}}.
$$

A {\it graded, enriched cycle $E^\bullet$\/} is simply an enriched cycle $E^i$ for $i$ in some bounded set of integers. An single
enriched cycle is considered as a graded enriched cycle by being placed totally in degree zero. The analytic set $V$ is a {\it
component of $E^\bullet$} if and only if $V$ is a component of $E^i$ for some $i$, and the underlying set of $E^\bullet$ is
$|E^\bullet|=\cup_i|E^i|$. If $R$ is a domain, then $E^\bullet$ yields an ordinary cycle
$[E^\bullet]^{\operatorname{ord}}:=\sum_i(-1)^i(\operatorname{rk}(E^i_V))[V]$. If $k$ is an integer, we define the {\it
$k$-shifted graded, enriched cycle}
$E^\bullet[k]$ by $(E^\bullet[k])^i:= E^{i+k}$.

If $q$ is a finitely-generated module and $E^\bullet$ is a graded enriched cycle, then we define the graded enriched cycle
$qE^\bullet$ by
$(qE^\bullet)^i:=\sum_V(q\otimes E^i_V)[V]$. The (direct) sum of two graded enriched cycles $D^\bullet$ and $E^\bullet$ is given
by $(D^\bullet+ E^\bullet)^i_V := D^i_V\oplus E^i_V$. If $D^i$ properly intersects $E^j$ for all $i$ and $j$, then we say that
$D^\bullet$ and $E^\bullet$ {\it intersect properly\/} and we define the intersection product by 
$$ (D^\bullet\odot E^\bullet)^k:=\sum_{i+j = k}(D^i\odot E^j).
$$

Whenever we use the enriched intersection product symbol, we mean that we are considering the objects on both sides of $\odot$ as graded, enriched cycles, even if we do not superscript by $\operatorname{enr}$ or  $\bullet$.
\end{defn}

\bigskip

Let $\tau: W\rightarrow Y$ be a proper morphism between analytic spaces. If $C = \sum n_V[V]$ is an ordinary positive cycle in
$W$, then the proper push-forward $\tau_*(C) = \sum n_V\tau_*([V])$ is a well-defined ordinary cycle. 

\begin{defn}\label{def:properpush} If $E^\bullet=\sum_V
E^\bullet_V[V]$ is an enriched cycle in $W$, then we define the {\bf proper push-forward of $E^\bullet$ by $\tau$} to be the
graded enriched cycle $\tau^\bullet_*(E^\bullet)$ defined by
$$
\tau^j_*(E^\bullet) \ := \ \sum_V E^j_V[\tau_*([V])]^{\operatorname{enr}}.
$$ 
\end{defn}

The ordinary projection formula for divisors ([{\bf F}], 2.3.c) immediately implies the following enriched version. 

\begin{prop}\label{prop:properpush} Let
$E^\bullet$ be a graded enriched cycle in $X$. Let $W:=|E^\bullet|$. Let $\tau:W\rightarrow Y$ be a proper morphism, and let
$g: Y\rightarrow\mathbb C$ be an analytic function such that
$g\circ\tau$ is not identically zero on any component of $E^\bullet$. Then, $g$ is not identically zero on any component of
$\tau^\bullet_*(E^\bullet)$ and
$$
\tau^\bullet_*\big(E^\bullet\ \odot \ V(g\circ\tau)\big) \ = \ \tau^\bullet_*(E^\bullet) \ \odot \ V(g).
$$
\end{prop}

\vskip .3in

\begin{defn}\label{def:gecc} Suppose that $\Fdot$ is a bounded complex of sheaves, which is constructible with respect to an
analytic Whitney stratification $\strat$, in which the strata are connected. For $S\in\strat$, let  $d_S:=\dim S$. If $(\mathbb
N_S, \mathbb L_S)$ is a pair consisting of a normal slice and complex link, respectively, to the stratum $S$,
then, for each integer $k$, the isomorphism-type of the module
$\mathbb H^{k-d_S}(\mathbb N_S, \mathbb L_S; \Fdot)$ is independent of the choice of $(\mathbb N_S, \mathbb L_S)$; we
refer to
$\mathbb H^{k-d_S}(\mathbb N_S,\mathbb L_S; \Fdot)$ as the {\bf degree $k$ Morse module of $S$ with respect to
$\Fdot$}. 

The {\bf graded, enriched characteristic cycle of $\Fdot$ in the cotangent bundle $T^*\U$} is defined in degree
$k$ to be 
$$
\gecc^k(\Fdot) :=\sum_{S\in\strat} H^{k-d_S}(\mathbb N_S,\mathbb L_S; \Fdot)\big[\,\overline{T^*_{{}_{S}}\mathcal U}\,\big].
$$
\end{defn}

\vskip .3in

\begin{rem} There are no canonical choices for defining the the normal slices or complex links of strata. However, as two enriched cycles are equal provided that the component modules are all isomorphic, the graded, enriched characteristic cycle is well-defined.
\end{rem}

\begin{exm}\label{exm:gecc} We wish a give a simple example of calculating a graded, enriched characteristic cycle.

Let $f:\C^3\rightarrow\C$ be given by $f(x, y, t) =y(y^2-x^3-t^2x^2)$, and let $X:=V(f)=V(y)\cup V(y^2-x^3-t^2x^2)$. The singular set of $X$, $\Sigma X$, is the $1$-dimensional set $V(x, y)\cup V(x+t^2, y)$.  Thus, near the origin (actually, in this specific example, globally), 
$$\strat:=\{V(y)- V(y^2-x^3-t^2x^2), V(y^2-x^3-t^2x^2)-V(y), V(x, y)-\{\0\},  V(x+t^2, y)-\{\0\}, \{\0\}\}$$
 is a Whitney stratification of $X$ with connected strata. Let $\Fdot:=\Z^\bullet_X[2]$ (we shall discuss the shift by $2$ below), which is constructible with respect to any Whitney stratification of $X$. We wish to calculate $\gecc^\bullet(\Fdot)$.

\smallskip

First, consider the $2$-dimensional strata. Let $S_1:=V(y)- V(y^2-x^3-t^2x^2)$. Then, $\N_{S_1}$ is simply a point, and $\cL_{S_1}$ is empty. Hence, $H^{k-2}(\N_{S_1}, \cL_{S_1}; \Fdot)=H^{k}(\N_{S_1}, \cL_{S_1}; \Z)$ isomorphic to $\Z$ if $k=0$, and is $0$ if $k\neq 0$. The same conclusion holds if $S_1$ is replaced by $S_2:=V(y^2-x^3-t^2x^2)-V(y)$.

Now, consider the $1$-dimensional strata. Let $S_3:= V(x, y)-\{\0\}$, and $S_4:=V(x+t^2, y)-\{\0\}$. The normal slice $\N_{S_3}$ is, as a germ, up to analytic isomorphism, three complex lines in $\C^2$, which intersect at a point, and $\cL_{S_3}$ is three points. Similarly, the normal slice $\N_{S_4}$ is, as a germ, up to analytic isomorphism, two complex lines in $\C^2$, which intersect at a point, and $\cL_{S_4}$ is two points. Hence, $H^{k-1}(\N_{S_3}, \cL_{S_3}; \Fdot)=H^{k+1}(\N_{S_3}, \cL_{S_3}; \Z)$ isomorphic to $\Z^2$ if $k=0$, and is $0$ if $k\neq 0$. Similarly, $H^{k-1}(\N_{S_4}, \cL_{S_4}; \Fdot)=H^{k+1}(\N_{S_4}, \cL_{S_4}; \Z)$ isomorphic to $\Z$ if $k=0$, and is $0$ if $k\neq 0$.

Finally, consider the stratum $\{\0\}$. Then, $\N_{\{\0\}}$ is all of $X$, intersected with a small ball around the origin. The complex link $\cL_{\{\0\}}$ is usually referred to as simply the complex link of $X$ at $\0$. Thus, $\cL_{\{\0\}}$ has the homotopy-type of a bouquet of $1$-spheres (see \cite{levan}), and the number of spheres in this bouquet is equal to the intersection number $(\Gamma^1_{f, L}\cdot V(L))_\0$, where $L$ is any linear form such that $d_\0L$ is not a degenerate covector from strata of $X$ at $0$ (see \cite{stratmorse}), and the relative polar curve here is the classical one from the beginning of \secref{sec:genle}. We claim that we may use $L:=t$ for this calculation.

To see this, first note that $V(y^2-x^3-t^2x^2)$ is the classic example of a space such that the regular part satisfies Whitney's condition (a) along the $t$-axis (or, alternatively, this is an easy exercise). Thus, $d_\0t$ is not a limit of conormals from $S_2$. Now, the closures of $S_1$, $S_3$, and $S_4$ are all smooth, and $d_\0t$ is not conormal to these closures at the origin.

To find the ordinary cycle $\Gamma^1_{f, t}$, we take the components of the cycle below which are {\bf not} contained in $\Sigma f$:
$$
V\left(\frac{\partial f}{\partial x}, \frac{\partial f}{\partial y}\right) = V(y(-3x^2-2t^2x), 3y^2-x^3-t^2x^2)=V(y, x^2(x+t^2)) + V(x(3x+2t^2), 3y^2-x^3-t^2x^2)=
$$
$$
2V(x, y)+V(x+t^2, y) +2V(x, y) + V(3x+2t^2, 3y^2-x^3-t^2x^2).
$$

\smallskip

\noindent Thus, $\Gamma^1_{f, t} = V(3x+2t^2, 3y^2-x^3-t^2x^2)$, and $(\Gamma^1_{f, t}\cdot V(t))_\0 = [V(3x+2t^2, 3y^2-x^3-t^2x^2, t)]_\0 = 2$,
and $H^{k-0}(\N_{\{\0\}}, \cL_{\{\0\}}; \Fdot)=H^{k+2}(\N_{\{\0\}}, \cL_{\{\0\}}; \Z)$ is isomorphic to $\Z^2$ if $k=0$, and is $0$ if $k\neq 0$.

\smallskip

Therefore, we find that $\gecc^k(\Fdot)= 0$ if $k\neq 0$, and 
$$
\gecc^0(\Fdot) = \Z\left[\overline{T^*_{S_1}\C^3}\right]+ \Z\left[\overline{T^*_{S_2}\C^3}\right] +\Z^2\left[\overline{T^*_{S_3}\C^3}\right]+ \Z\left[\overline{T^*_{S_4}\C^3}\right]+\Z^2 \left[T^*_{\{\0\}}\C^3\right].
$$

The fact that $\gecc^\bullet(\Fdot)$ is concentrated in degree $0$ is equivalent to the fact that $\Z^\bullet_X[2]$ is a perverse sheaf (see \cite{singenrich}), and was the reason for including the shift by $2$. The constant sheaf on any connected, local complete intersection, shifted by the dimension of the space, is perverse.

\bigskip

The reader is invited to take the most simple space $Y$ which is not a local complete intersection -- two planes $P_1$ and $P_2$ in $\C^4$, which intersect at only the origin -- and show that, if $\Adot=\Z^\bullet_Y[2]$, then
$$
\gecc^0(\Adot) = \Z\big[T^*_{P_1}\C^4\big]+\Z\big[T^*_{P_2}\C^4\big],
$$
$$
\gecc^{-1}(\Adot) = \Z\big[T^*_{\{\0\}}\C^4\big],
$$
and $\gecc^k(\Adot)=0$ for $k\neq 0, -1$.
\end{exm}

\section{The Main Definitions and First Main Theorem}\label{sec:main}

Throughout the remainder of this paper, we will use the notation established in \secref{sec:genle}:  $\U$ is an open neighborhood of the origin in $\C^{n+1}$, $\tilde f$ and $\tilde g$ are analytic functions from $(\U, \0)$ to $(\C, 0)$, $X$ is a complex analytic subset of $\U$, $f$ and $g$ denote the restrictions of $\tilde f$ and $\tilde g$, respectively, to $X$, and $\strat$ is a Whitney stratification of $X$, with connected strata, such that $V(f)$ is a union of strata. 

We use $(z_0, \dots, z_n)$ for coordinates on $\U$, and identify $T^*\U$ with $\U\times\C^{n+1}$, using $(w_0, \dots, w_n)$ for cotangent coordinates, so that $(p, w_0d_pz_0+\dots+w_nd_pz_n)$ corresponds to $(p, (w_0, \dots, w_n))$. Let $\pi:T^*\U\rightarrow\U$ denote the projection.  Below, we consider the image, $\im d\tilde g$, of $d\tilde g$ in $T^*\U$; this scheme is defined by 
$$V\left(w_0- \frac{\partial \tilde g}{\partial z_0}, \dots, w_n- \frac{\partial \tilde g}{\partial z_n}\right)\subseteq \U\times\C^{n+1}.$$
We will consider $\im d\tilde g$ as a scheme, an analytic set, an ordinary cycle, and as a graded, enriched cycle; we will denote all of these by simply $\im d\tilde g$, and explicitly state what structure we are using or let the context make the structure clear.

We do not require our base ring to be $\Z$ (as we did in \secref{sec:genle}). We let $R$, our base ring, be any regular, Noetherian ring with finite Krull dimension (e.g., $\Z$, $\Q$, or $\C$). This implies that every finitely-generated $R$-module has finite projective dimension (in fact, it implies that the projective dimension of the module is at most $\dim R$). We let $\Fdot$ be a bounded, constructible complex of sheaves of $R$-modules on $X$. Let $\strat(\Fdot):=\{S\in\strat\ |\ \hyp^*(\mathbb N_S,\mathbb L_S; \Fdot)\neq 0\}$; we refer to the elements of $\strat(\Fdot)$ as the {\it $\Fdot$-visible strata of $\strat$}.

\bigskip

Suppose that $M$ is a complex submanifold of $\U$. Recall:

\begin{defn}\label{def:main} The relative conormal space $T^*_{\tilde f_{|_M}}\U$ is given by 
$$
T^*_{\tilde f_{|_M}}\U :=\{(x, \eta)\in T^*\U\ |\ \eta(T_xM\cap \ker d_x\tilde f)=0\}.
$$

If $M\subseteq X$, then $T^*_{\tilde f_{|_M}}\U$ depends on $f$, but not on the particular extension $\tilde f$. In this case, we write $T^*_{f_{|_M}}\U$ in place of $T^*_{\tilde f_{|_M}}\U$.
\end{defn}

\begin{defn}\label{def:relconorm} The {\bf graded, enriched relative conormal cycle, $\big(T^*_{{}_{f,\Fdot}}\U\big)^\bullet$, of $f$, with respect to $\Fdot$}, is defined by
$$\big(T^*_{{}_{f,\Fdot}}\U\big)^k:=\sum_{\substack{S\in\strat(\Fdot)\\ f_{|_S}\neq{\rm\ const.}}}H^{k-d_S}(\mathbb N_S,\mathbb L_S; \Fdot)\left[\overline{T^*_{f_{|_S}}\U}\right].$$
\end{defn}

\begin{exm}\label{exm:gercc} Let us return to the setting of \exref{exm:gecc}, where $X=V(y)\cup V(y^2-x^3-t^2x^2)$ and $\Fdot=\Z^\bullet_X[2]$. We had Whitney strata consisting of $\{\0\}$, $S_1=V(y)- V(y^2-x^3-t^2x^2)$, $S_2=V(y^2-x^3-t^2x^2)-V(y)$, $S_3= V(x, y)-\{\0\}$, and $S_4=V(x+t^2, y)-\{\0\}$.

We found that
$\gecc^k(\Fdot)= 0$ if $k\neq 0$, and 
$$
\gecc^0(\Fdot) = \Z\left[\overline{T^*_{S_1}\C^3}\right]+ \Z\left[\overline{T^*_{S_2}\C^3}\right] +\Z^2\left[\overline{T^*_{S_3}\C^3}\right]+ \Z\left[\overline{T^*_{S_4}\C^3}\right]+\Z^2 \left[T^*_{\{\0\}}\C^3\right].
$$
We will calculate $\big(T^*_{{}_{x,\Fdot}}\C^3\big)^\bullet$.

\smallskip

As we said above, we identify $T^*\C^3$ with $\C^3\times\C^3$, and will use coordinates $(w_0, w_1, w_2)$ for cotangent coordinates, so that $(w_0, w_1, w_2)$ represents $w_0dx+w_1dy+w_2dt$.

Since $x$ is identically zero on $\{\0\}$ and $S_3$, these two strata are not used in the calculation of $\big(T^*_{{}_{x,\Fdot}}\C^3\big)^\bullet$. For the $1$-dimensional stratum $S_4$, $\left[\overline{T^*_{x_{|_{S_4}}}\C^3}\right]$ is the $4$-dimensional cycle $V(x+t^2,y)\subseteq \C^3\times\C^3$.

The fiber of $T^*_{x_{|_{S_1}}}\C^3$ over any $p\in S_1$ is  
$$(T^*_{S_1}\C^3)_p+<d_px>:=\{\omega+ad_px\ |\ \omega\in (T^*_{S_1}\C^3)_p, a\in\C\}=\{bd_py+ad_px\ |\ a,b\in\C\}.$$ 
Hence, $\left[\overline{T^*_{x_{|_{S_1}}}\C^3}\right]= V(y, w_2)$.

The fiber of $T^*_{x_{|_{S_2}}}\C^3$ over any $p\in S_2$ which is a regular point of $x$ restricted to $S_2$ is
$$(T^*_{S_2}\C^3)_p+<d_px>:=$$
$$\{\omega+ad_px\ |\ \omega\in (T^*_{S_2}\C^3)_p, a\in\C\}=\{b\big((-3x^2-2t^2x)d_px+2yd_py-2tx^2d_pt\big)+ad_px\ |\ a,b\in\C\}.$$
The form $w_0d_px+w_1d_py+w_2d_pt$ is in this set if and only if the determinant of the following matrix is $0$:
$$\left[\begin{matrix}w_0 & w_1 & w_2\\ -3x^2-2t^2x & 2y & -2tx^2\\ 1 & 0 & 0
\end{matrix}\right],$$
i.e., if and only if $yw_2+tx^2w_1=0$.
It is tempting to conclude that $\left[\overline{T^*_{x_{|_{S_2}}}\C^3}\right]$ equals $V(y^2-x^3-t^2x^2, yw_2+tx^2w_1)$, but this is not the case; we must eliminate any components of $V(y^2-x^3-t^2x^2, yw_2+tx^2w_1)$ which are contained in $V(x,y)$. Our notation for the resulting scheme (a gap sheaf, see \cite{numcontrol}, I.1) is 
$$V(y^2-x^3-t^2x^2, yw_2+tx^2w_1)\lnot V(x,y).$$
Note that, as schemes,
$$V(y^2-x^3-t^2x^2, yw_2+tx^2w_1)= V(y^2-x^3-t^2x^2, yw_2+tx^2w_1, y^2w_2+ytx^2w_1)=$$
$$
V(y^2-x^3-t^2x^2, yw_2+tx^2w_1, (x^3+t^2x^2)w_2+ytx^2w_1).
$$
Using \cite{numcontrol}, I.1.3.iv, we find that, as cycles,
$$V(y^2-x^3-t^2x^2, yw_2+tx^2w_1)\lnot V(x,y) = V(y^2-x^3-t^2x^2, yw_2+tx^2w_1, (x+t^2)w_2+ytw_1).$$
(This last equality need not be true on the level of schemes, since our generators do not form a regular sequence and, hence, there may be embedded subvarieties.)

Therefore, we find that $\big(T^*_{{}_{x,\Fdot}}\C^3\big)^k$ is $0$ unless $k=0$, and 
$$
\big(T^*_{{}_{x,\Fdot}}\C^3\big)^0=\Z[V(y, w_2)]\ +\ \Z[V(y^2-x^3-t^2x^2, yw_2+tx^2w_1, (x+t^2)w_2+ytw_1)]\ +\ \Z[V(x+t^2, y)].
$$
\end{exm}

\bigskip

We now wish to define the graded, enriched relative polar curve.  Note that the projection $\pi$ induces an isomorphism from the analytic set $\im d\tilde g$ to $\U$. We will use the proper push-forward (\defref{def:properpush}) of the map $\pi$ restricted to $\im d\tilde g$; we will continue to denote this restriction by simply $\pi$. 

  By our conventions in \secref{sec:enrcycle}, the graded, enriched $\im d\tilde g$ is zero outside of degree $0$, and is $R[\im d\tilde g]$ in degree $0$.

\begin{defn}\label{def:polarcurve} If $S\in\strat$ and $f_{|_S}$ is not constant, we define the  {\bf relative polar set}, $\big|\Gamma_{f,\tilde g}(S)\big|$, to be $\pi\left(\overline{T^*_{f_{|_S}}\U}\ \cap\ \im d\tilde g\right)$; if this set is purely $1$-dimensional, so that $\overline{T^*_{f_{|_S}}\U}$ and $\im d\tilde g$ intersect properly, we define the (ordinary) {\bf relative polar cycle}, $\Gamma_{f,\tilde g}(S)$, to be the cycle $\pi_*\left(\left[\overline{T^*_{f_{|_S}}\U}\right]\cdot [\im d\tilde g]\right)$.

The {\bf relative polar set}, $\big|\Gamma_{f, \tilde g}(\Fdot)\big|$, is defined by 
$$
\big|\Gamma_{f, \tilde g}(\Fdot)\big|:= \pi\left(\big|\big(T^*_{{}_{f,\Fdot}}\U\big)^\bullet\big|\cap \im d\tilde g\right).
$$

Each $1$-dimensional component $C$ of $\big|\Gamma_{f, \tilde g}(\Fdot)\big|$ is the image of a component of $\big|\big(T^*_{{}_{f,\Fdot}}\U\big)^\bullet\big|\cap \im d\tilde g$ along which $\big|\big(T^*_{{}_{f,\Fdot}}\U\big)^\bullet\big|$ and  $\im d\tilde g$ intersect properly. We give such a component $C$ the structure of the graded, enriched cycle whose underlying set is $C$ and whose graded, enriched cycle structure is given by $\pi_*^\bullet\left(\big(T^*_{{}_{f,\Fdot}}\U\big)^\bullet\odot \im d\tilde g\right)$ over generic points in $C$. We refer to this as the {\bf graded, enriched cycle structure of $C$ in $\big|\Gamma_{f, \tilde g}(\Fdot)\big|$}.

If  $\big|\Gamma_{f, \tilde g}(\Fdot)\big|$ is purely $1$-dimensional, we say that the {\bf graded, enriched relative polar curve}, $\big(\Gamma^1_{f, \tilde g}(\Fdot)\big)^\bullet$, is defined, and is given by 
$$
\big(\Gamma^1_{f, \tilde g}(\Fdot)\big)^\bullet:= \pi_*^\bullet\left(\big(T^*_{{}_{f,\Fdot}}\U\big)^\bullet\odot \im d\tilde g\right),
$$
i.e., 
$$
\big(\Gamma^1_{f, \tilde g}(\Fdot)\big)^k = \sum_{\substack{S\in\strat(\Fdot)\\ f_{|_S}\neq{\rm\ const.}}}H^{k-d_S}(\mathbb N_S,\mathbb L_S; \Fdot)\left(\Gamma_{f, \tilde g}(S)\right)^{\operatorname{enr}}.
$$
\end{defn}

\begin{rem}\label{rem:polar} In the notation for the polar curve, we write $\tilde g$, not simply $g$; we do not, in fact, know if $\big(\Gamma^1_{f, \tilde g}(\Fdot)\big)^\bullet$ is independent of the extension to $\tilde g$. However, \thmref{thm:main1} will imply that, when $\big(\Gamma^1_{f, \tilde g}(\Fdot)\big)^\bullet$ is defined and has no component on which $f$ is constant, then $\big(\Gamma^1_{f, \tilde g}(\Fdot)\big)^\bullet$ is independent of the extension $\tilde g$. It is also not difficult to show that the set $\big|\Gamma_{f, \tilde g}(\Fdot)\big|$ is independent of the extension of $g$, but we shall not need this result here.

Note that $\overline{T^*_{f_{|_S}}\U}\cap \im d\tilde g$ is at least $1$-dimensional at each point of intersection, and so $\big|\Gamma_{f, \tilde g}(\Fdot)\big|$ has no isolated points. Also, note that, as $\big|\big(T^*_{{}_{f,\Fdot}}\U\big)^\bullet\big|\cap \im d\tilde g$ is a closed subset of $\im d\tilde g$, and $\pi$ induces an isomorphism from $\im d\tilde g$ to $\U$, $\big|\Gamma_{f, \tilde g}(\Fdot)\big|$ is a closed subset of $\U$.

Finally, the reader may wonder about the symmetry of our definition. It is {\bf not} true for arbitrary $\tilde f$ and $\tilde g$ that even the sets  $\big|\Gamma_{f, \tilde g}(\Fdot)\big|$ and  $\big|\Gamma_{g, \tilde f}(\Fdot)\big|$ are equal; see \remref{rem:symm2}. However, \propref{prop:sigma} will imply that the components of these two sets along which neither $f$ nor $g$ are constant are the same. Hence, we refer to a component of $\big|\Gamma_{f, \tilde g}(\Fdot)\big|$ along which neither $f$ nor $g$ is is constant as a {\it symmetric component of $\big|\Gamma_{f, \tilde g}(\Fdot)\big|$}.

By moving to a generic point $p$ on a $1$-dimensional symmetric component $C$ of $\big|\Gamma_{f, \tilde g}(\Fdot)\big|$ and applying \corref{cor:main2}, one can show that the graded, enriched cycle structure of $C$ in $\big|\Gamma_{f, \tilde g}(\Fdot)\big|$ is the same as that of $C$ in $\big|\Gamma_{g, \tilde f}(\Fdot)\big|$.
\end{rem}

\medskip

\begin{exm}\label{ex:gerpc} We continue with our setting from \exref{exm:gecc} and \exref{exm:gercc}, and consider $X=V(y)\cup V(y^2-x^3-t^2x^2)$ and $\Fdot=\Z^\bullet_X[2]$. We will calculate $\big(\Gamma^1_{x, t}(\Fdot)\big)^\bullet$.

Using the isomorphism $T^*\C^3\cong \C^3\times\C^3$ from \exref{exm:gercc}, $\im dt$ is the scheme 
$$V\left(w_0-\frac{\partial t}{\partial x}, \ w_1-\frac{\partial t}{\partial y}, \ w_2-\frac{\partial t}{\partial t} \right)=V(w_0, w_1, w_2-1).$$

In \exref{exm:gercc}, we found that $\big(T^*_{{}_{x,\Fdot}}\C^3\big)^k$ is $0$ unless $k=0$, and 
$$
\big(T^*_{{}_{x,\Fdot}}\C^3\big)^0=\Z[V(y, w_2)]\ +\ \Z[V(y^2-x^3-t^2x^2, yw_2+tx^2w_1, (x+t^2)w_2+ytw_1)]\ +\ \Z[V(x+t^2, y)].
$$
Let us write $E$ for the cycle $V(y^2-x^3-t^2x^2, yw_2+tx^2w_1, (x+t^2)w_2+ytw_1)$ throughout the remainder of this example.

Thus, $\big(\Gamma^1_{x, t}(\Fdot)\big)^k$ is $0$ unless $k=0$ and, to calculate $\big(\Gamma^1_{x, t}(\Fdot)\big)^0$, we need first to calculate the three ordinary cycles
$$\pi_*\big(V(y, w_2)\cdot V(w_0, w_1, w_2-1)\big),$$
$$\pi_*\big(E\cdot V(w_0, w_1, w_2-1)\big),$$
and
$$\pi_*\big(V(x+t^2, y)\cdot V(w_0, w_1, w_2-1)\big).$$

Now, $V(y, w_2)\cap V(w_0, w_1, w_2-1)=\emptyset$, and so $\pi_*\big(V(y, w_2)\cdot V(w_0, w_1, w_2-1)\big)=0$. In addition, it is trivial that there is an equality of cycles $\pi_*\big(V(x+t^2, y)\cdot V(w_0, w_1, w_2-1)\big)= V(x+t^2, y)$. However, the remaining cycle is more difficult to calculate.

The difficulty in calculating 
$$\pi_*\big(E\cdot V(w_0, w_1, w_2-1)\big)$$
is due to the fact that $y^2-x^3-t^2x^2, yw_2+tx^2w_1, (x+t^2)w_2+ytw_1$ is not a regular sequence. To ``fix'' this, note that, in \exref{exm:gercc}, we saw that, as cycles, there is an equality
$$
V(y^2-x^3-t^2x^2, yw_2+tx^2w_1)= C+E,
$$
where the underlying set $|C|\subseteq V(x, y)$. Now, it is trivial that, as sets, 
$$E\cap V(w_0, w_1, w_2-1)= V(x+t^2, y, w_0, w_1, w_2-1).
$$
Therefore,
$$
C\cdot V(w_0, w_1, w_2-1) + E\cdot V(w_0, w_1, w_2-1) = V(y^2-x^3-t^2x^2, yw_2+tx^2w_1)\cdot V(w_0, w_1, w_2-1) = 
$$
$$
V(y^2-x^3-t^2x^2, yw_2+tx^2w_1, w_0, w_1, w_2-1)= V(x^2(x+t^2), y, w_0, w_1, w_2-1)
$$
$$
2V(x,y,w_0, w_1, w_2-1)+ V(x+t^2, y, w_0, w_1, w_2-1).
$$
Thus, as cycles,
$$
E\cdot V(w_0, w_1, w_2-1) = V(x+t^2, y, w_0, w_1, w_2-1),
$$
and so $\pi_*(E\cdot V(w_0, w_1, w_2-1))= V(x+t^2, y)$.

Finally, we find that
$$
\big(\Gamma^1_{x, t}(\Fdot)\big)^0= \pi_*^0\left((T^*_{{}_{f,\Fdot}}\U\big)^\bullet\odot \im dt\right) = \Z[V(x+t^2, y)]+  \Z[V(x+t^2, y)] =  \Z^2[V(x+t^2, y)].
$$
\end{exm}

\bigskip

Before we can prove our main theorem of this section, we must recall three results from \cite{singenrich}.

\begin{thm}\label{thm:psigecc} {\rm (\cite{singenrich}, Theorem 3.3)} There is an equality of graded enriched cycles given by
$$
\gecc^\bullet(\psi_f[-1]\Fdot) \ = \ \big(T^*_{{}_{f, \Fdot}}\U\big)^\bullet\odot(V(f)\times\C^{n+1}).
$$ 
\end{thm}

\smallskip

We state the next two theorems for complexes of sheaves on $V(f)$, since that is the case in which we shall use them.

\begin{thm}\label{thm:isovan}{\rm (\cite{singenrich}, Theorem 3.4)} Let $\Adot$ be a bounded, constructible complex of sheaves of $R$-modules on $V(f)$. Then, $\dim_\0 (\supp\phi_{g}[-1]\Adot)\leq 0$ if and only if
$\dim_{\0}\pi\big(|\gecc^\bullet(\Adot)|\cap \operatorname{im}d\tilde g\big)\leq 0$, and when this
is the case, $\dim_{(\0, d_\0\tilde g)}\big(|\gecc^\bullet(\Adot)|\cap \operatorname{im}d\tilde g\big)\leq 0$ and 
$$H^k(\phi_{g}[-1]\Adot)_\0\cong \big(\gecc^k(\Adot)\ \odot\ \operatorname{im}d\tilde g\big)_{(\0, d_\0\tilde g)}.
$$
\end{thm}

\begin{thm}\label{thm:vansupp} {\rm (\cite{singenrich}, Theorem 3.5)} Let $\Adot$ be a bounded, constructible complex of sheaves of $R$-modules on $V(f)$. Then, there is an equality of sets given by
$$\bigcup_{v\in\C}\supp\phi_{g-v}[-1]\Adot\ =\pi\big(|\gecc^\bullet(\Adot)|\cap\operatorname{im}d\tilde g\big).
$$ 
\end{thm}

\bigskip

We need a lemma before we prove our first main theorem.

\bigskip

\begin{lem}\label{lem:phipsi} There is an equality of sets
$$
\supp\phi_g[-1]\psi_f[-1]\Fdot = V(f,g)\cap\big|\Gamma_{f, \tilde g}(\Fdot)\big|,
$$
and, in a neighborhood of the origin,
$V(f)\cap\big|\Gamma_{f, \tilde g}(\Fdot)\big|\subseteq V(g)$.
\end{lem}
\begin{proof} The equality follows from \thmref{thm:vansupp} by letting $\Adot=\psi_f[-1]\Fdot$, applying \thmref{thm:psigecc}, and then intersecting $V(g)$ with both sides of the equation from  \thmref{thm:vansupp}. The containment also follows from \thmref{thm:vansupp} by letting $\Adot=\psi_f[-1]\Fdot$, applying \thmref{thm:psigecc}, and then using that, near a point $p$ where $g(p)=0$, $\supp\phi_{g-v}[-1]\Adot$ is empty if $v\neq 0$.
\end{proof}

\smallskip

\begin{rem}\label{rem:prepolar} In much of our work, particularly in \cite{lecycles} and \cite{numcontrol}, we have used the notion of a {\it prepolar slice}. Using our current terminology, the condition that $V(\tilde g)$ is a prepolar slice for $f$ at $\0$ would be replaced with $\dim_\0 V(g)\cap \big|\Gamma_{f, \tilde g}(\Fdot)\big|\leq 0$. Note that, by \lemref{lem:phipsi}, this implies that $\dim_\0 V(f)\cap \big|\Gamma_{f, \tilde g}(\Fdot)\big|\leq 0$.
\end{rem}

\bigskip

We now prove our first main theorem.

\begin{thm}\label{thm:main1} The following are equivalent:
\begin{enumerate}
\item $\dim_\0\supp \phi_g[-1]\psi_f[-1]\Fdot\leq 0$;
\item $\dim_\0 V(f)\cap \big|\Gamma_{f, \tilde g}(\Fdot)\big|\leq 0$;
\item $\dim_\0 V(f, g)\cap \big|\Gamma_{f, \tilde g}(\Fdot)\big|\leq 0$;
\end{enumerate}
and, when these equivalent conditions hold, $\big(\Gamma^1_{f, \tilde g}(\Fdot)\big)^\bullet$ exists and 
$$H^k(\phi_g[-1]\psi_f[-1]\Fdot)_\0\cong \big(\big(\Gamma^1_{f, \tilde g}(\Fdot)\big)^k\odot V(f)\big)_\0,
$$
i.e., 
$$H^{k-1}(B_\epsilon\cap f^{-1}(\xi), B_\epsilon\cap f^{-1}(\xi)\cap g^{-1}(\nu);\ \Fdot)\ \cong\ \bigoplus_{S\in \strat^o(\Fdot)} \big(\mathbb H^{k-d_S }(\N_S, \cL_S  ; \Fdot)\big)^{j_S} , $$ 
where $0<|\xi|\ll |\nu|<\ll\epsilon\ll 1$, and 
$j_S = \big(\Gamma^1_{{}_{f, \tilde g}}(S)\cdot V(f)\big)_\mathbf 0$.
\end{thm}
\begin{proof} The equivalence of the conditions follows immediately from the lemma. Assume now that these conditions hold.

By \thmref{thm:isovan},
$$H^k(\phi_{g}[-1]\psi_f[-1]\Fdot)_\0\cong \big(\gecc^k(\psi_f[-1]\Fdot)\ \odot\ \operatorname{im}d\tilde g\big)_{(\0, d_\0\tilde g)}.
$$
Applying \thmref{thm:psigecc}, we find that
$$H^k(\phi_g[-1]\psi_f[-1]\Fdot)_{\0} \cong \big(\big(T^*_{{}_{f, \Fdot}}\U\big)^k\odot (V(f)\times\C^{n+1})\ \odot\ \im d\tilde g\big)_{(\0, d_\0\tilde g)}\cong$$
$$
\left(\pi^k_*\Big(\big(T^*_{{}_{f, \Fdot}}\U\big)^\bullet\odot (V(f)\times\C^{n+1})\ \odot\ \im d\tilde g\Big)\right)_\0,
$$
where this last isomorphism follows from the definition of the proper push-forward. By \propref{prop:properpush}, this last quantity equals
$$
\big(\pi^k_*\big(\big(T^*_{{}_{f, \Fdot}}\U\big)^\bullet\ \odot\ \im d\tilde g\big)\odot V(f)\big)_\0,
$$
which, by definition of the graded, enriched relative polar curve is equal to $\big(\big(\Gamma^1_{f, \tilde g}(\Fdot)\big)^k\odot V(f)\big)_\0$.
\end{proof}

\bigskip

We would like to know, of course, that the equivalent hypotheses of \thmref{thm:main1} are satisfied in the classical case where $f$ is fixed and $\tilde g$ is chosen to be a generic linear form.

\smallskip

\begin{prop}\label{prop:geng} 
\begin{enumerate}
\item There exists a non-zero linear form $\mathfrak l$ such that $\0\not\in\big|\Gamma_{f, \mathfrak l}(\Fdot)\big|$ if and only if for generic linear $\mathfrak l$, $\0\not\in\big|\Gamma_{f, \mathfrak l}(\Fdot)\big|$.

\item For generic linear $\mathfrak l$, $\dim_\0 V(f)\cap \big|\Gamma_{f, \mathfrak l}(\Fdot)\big|\leq 0$ and $\dim_\0 V(\mathfrak l)\cap \big|\Gamma_{f, \mathfrak l}(\Fdot)\big|\leq 0$.
\end{enumerate}
\end{prop}
\begin{proof} The proof of Item 1 is standard. Suppose that there exists a non-zero linear form $\mathfrak l$ such that $\0\not\in\big|\Gamma_{f, \mathfrak l}(\Fdot)\big|$. Then, the projective class $[d_\0\mathfrak l]$ is not in the fiber $\Proj\big(\big|\big(T^*_{{}_{f,\Fdot}}\U\big)^\bullet\big|\big)_\0\subseteq\Proj^n$. Thus, $\Proj\big(\big|\big(T^*_{{}_{f,\Fdot}}\U\big)^\bullet\big|\big)_\0$ is a proper analytic subset of $\Proj^n$. This implies Item 1.

\smallskip

\noindent Proof of Item 2:

If $\Adot$ is any bounded, constructible complex of sheaves on any complex analytic $Y\subset\U$ and $\0\in Y$, then for generic linear $\mathfrak l$ on $U$, $\0$ is an isolated  point in $\supp\phi_l[-1]\Adot$. This is well-known; see, for instance, Theorem 2.4  of \cite{pervcohovan}. Thus, for generic $\mathfrak l$, $\dim_\0\supp \phi_{\mathfrak l}[-1]\psi_f[-1]\Fdot\leq 0$, which, by \thmref{thm:main1}, is equivalent to $\dim_\0 V(f)\cap \big|\Gamma_{f, \mathfrak l}(\Fdot)\big|\leq 0$.

Now, refine, if necessary, our Whitney stratification $\strat$ to obtain a Whitney stratification $\strat^\prime$ which also satisfies Thom's $a_f$ condition. For generic linear $\mathfrak l$, $V(\mathfrak l)$ will transversely intersect all of the strata of $\strat^\prime$ in a neighborhood of the origin, except possibly at the origin itself. Fix such an $\mathfrak l$. We claim that $\dim_\0 V(\mathfrak l)\cap \big|\Gamma_{f, \mathfrak l}(\Fdot)\big|\leq 0$.

Since $\strat^\prime$ is an $a_f$ stratification, $
\bigcup_{S\in\strat^\prime}\overline{T^*_{f_{|_S}}\U} = \bigcup_{S\in\strat^\prime}T^*_{f_{|_S}}\U$. Our choice of $\mathfrak l$ implies that, if $p\neq \0$ and $p\in S\in\strat^\prime$, then $d_p\mathfrak l\not\in (T^*_S\U)_p$. Suppose now that we have an analytic path $\alpha(t)$ such that $\alpha(0)=\0$ and, for $t\neq 0$, $\alpha(t)\neq\0$ and $\alpha(t)\in V(\mathfrak l)\cap \big|\Gamma_{f, \mathfrak l}(\Fdot)\big|$. We wish to arrive at a contradiction.

If $f(\alpha(t))=0$ when $|t|$ is small, then we are finished, since $\dim_\0 V(f)\cap \big|\Gamma_{f, \mathfrak l}(\Fdot)\big|\leq 0$. So, assume that for $|t|$ small and non-zero, $f(\alpha(t))\neq 0$.

For $|t|$ small and non-zero, $\alpha(t)$ is contained in a single stratum $S\in\strat^\prime$. Near the origin, the $\strat^\prime$-stratified critical locus is contained in $V(f)$; hence, by the assumption in the previous paragraph, for $|t|$ small and non-zero, $d_{\alpha(t)}\tilde f\not\in (T^*_S\U)_{\alpha(t)}$. From our definition of $\alpha(t)$, and the discussion two paragraphs above, it follows that, for $|t|$ small and non-zero, $d_{\alpha(t)}\mathfrak l \not\in (T^*_S\U)_{\alpha(t)}$, $\mathfrak l(\alpha(t))\equiv 0$, and $d_{\alpha(t)}\mathfrak l\in \big(T^*_{f_{|_S}}\U\big)_{\alpha(t)}= (T^*_S\U)_{\alpha(t)}+<d_{\alpha(t)}\tilde f>$, where the last equality uses that $d_{\alpha(t)}\tilde f\not\in (T^*_S\U)_{\alpha(t)}$.

 Thus, for $|t|$ small and non-zero, there exists $c_t\in \C$ such that 
$$
d_{\alpha(t)}\mathfrak l+c_td_{\alpha(t)}\tilde f\in (T^*_S\U)_{\alpha(t)}\leqno(\dagger).
$$
Evaluating at $\alpha^\prime(t)$, and using that $\mathfrak l(\alpha(t))\equiv 0$ and $\alpha^\prime(t)\in T_{\alpha(t)}S$, we immediately conclude that $c_t\big(f(\alpha(t))\big)^\prime\equiv 0$. However, $c_t$ cannot be zero, for otherwise $(\dagger)$ would imply that $d_{\alpha(t)}\mathfrak l\in (T^*_S\U)_{\alpha(t)}$. Therefore, we must have that $\big(f(\alpha(t))\big)^\prime\equiv 0$, which implies that $f(\alpha(t))\equiv 0$, since $f(\alpha(0))=0$. This is a contradiction.
\end{proof}

\bigskip

\begin{exm} We continue where we left off in \exref{ex:gerpc}: $X=V(y)\cup V(y^2-x^3-t^2x^2)$ and $\Fdot=\Z^\bullet_X[2]$. We found that $\big(\Gamma^1_{x, t}(\Fdot)\big)^\bullet$ was concentrated in degree $0$, and
$$
\big(\Gamma^1_{x, t}(\Fdot)\big)^0=   \Z^2[V(x+t^2, y)].
$$
Thus, \thmref{thm:main1} tells us that $H^k(\phi_t[-1]\psi_x[-1]\Z^\bullet_X[2])_\0$ is $0$ unless $k=0$, and
$$
H^0(\phi_t[-1]\psi_x[-1]\Z^\bullet_X[2])_\0\cong \left( \Z^2[V(x+t^2, y)]\odot V(x)\right)_\0=\Z^4.
$$
\end{exm}

\section{The Discriminant as a Complex of Sheaves}\label{sec:discrim}

\bigskip

\thmref{thm:main1}, and its elegant, formal proof, was our motivation for defining the graded, enriched relative polar curve as we did. Of course, it would be nice to have a generalization of the result of L\^e in its original form, as it appears in \thmref{thm:old}: a result which gives $\hyp^*(F_{f, \mathbf 0}, F_{f_{|_{V(g)}}, \mathbf 0} ; \Fdot)$. In fact, we could easily prove such a result by appealing to the discriminant and Cerf diagram, if only we could push the complex $\Fdot$ down to the discriminant in some nice way.

There is one serious technical issue involved: we must show that a suitable neighborhood of origin pushes down by $(g, f)$ to a {\bf complex} analytically constructible complex, a complex which is constructible with respect to a stratification which is essentially determined by the image of the enriched relative polar curve. The main problem is that, on an open neighborhood of the origin, $(g, f)$ will not be a proper map and, if we instead use a domain with boundary on which $(g, f)$ is proper, then the boundary causes us to leave the complex analytic setting.

This is precisely the type of problem that is addressed by the microlocal theory of Kashiwara and Schapira in \cite{kashsch}, and we will use the micro-support of complexes of sheaves on real semianalytic sets. It will take a fair amount of preliminary work before we arrive at the desired result.

\bigskip

Suppose that $M$ is a $C^\infty$ manifold, $Z$ is a subspace of $M$, and $\Adot$ is a bounded complex of sheaves of $R$-modules on $Z$. Then, Kashiwara and Schapira define the micro-support, $SS(\Adot)\subseteq T^*M$, of $\Adot$ in 5.1.2 of \cite{kashsch}. Intuitively, $(p, \eta)\in SS(\Adot)$ if and only if the local hypercohomology of $Z$, with coefficients in $\Adot$, changes as one ``moves'' in the direction of $\eta$. 

In our fixed complex analytic setting, where $X$ is a complex analytic subset of $\U$ and $\Fdot$ is complex analytically constructible, the micro-support is easy to describe.

\begin{prop} {\rm (\cite{micromorse}, Theorem 4.13)} 
$$SS(\Fdot)=|\gecc^\bullet(\Fdot)|=\bigcup_{S\in\strat(\Fdot)}\overline{T^*_S\U}.$$
\end{prop}

\smallskip

We need to define the critical locus of complex analytic maps relative to the complex $\Fdot$. If $p\in X$, we shall write $SS_p(\Fdot)$ for the fiber $\pi^{-1}(p)\cap SS(\Fdot)$.

\begin{defn}\label{def:critlocus} The {\bf $\Fdot$-critical locus of $f$}, $\Sigma_{\Fdot}f$, is the set $\{p\in X\ |\ H^*(\phi_{f-f(p)}[-1]\Fdot)_p\neq 0\}$.
\end{defn}

\begin{prop}\label{prop:fdotcrit}{\rm(Theorem 2.4 and Remark 2.5 of \cite{pervcohovan})} The closure $\overline{\Sigma_{\Fdot}f}$ is equal to $\pi\big(\im d\tilde f\cap SS(\Fdot)\big)$ and $f$ is constant along the components of  $\overline{\Sigma_{\Fdot}f}$.
\end{prop}

\medskip

We need to generalize $\overline{\Sigma_{\Fdot}f}$ to the case where $f$ is a real analytic map whose codomain has dimension greater than one,  and where we replace $\Fdot$ by something more general.

\medskip

We may consider $T^*\U$ with its complex analytic structure, as we have been up to this point, or with its real analytic structure. When it is important for us to distinguish these structures, we will write  $(T^*\U)^\C$ and $(T^*\U)^\R$, respectively, and we remind the reader that, for $p\in\U$, there is an $\R$-linear isomorphism from $(T^*\U)^\C_p$ to $(T^*\U)^\R_p$ given by mapping $\eta$ to the real part $\operatorname{Re}\eta$ (or $2\operatorname{Re}\eta$). If $\eta_1, \dots, \eta_k\in (T^*\U)^\C_p$, this isomorphism identifies the complex span $\langle \eta_1, \dots, \eta_k\rangle^\C$ with the real span $\langle \operatorname{Re}\eta_1,  \operatorname{Im}\eta_1, \dots, \operatorname{Re}\eta_k,  \operatorname{Im}\eta_k\rangle^\R$. When the structure is clear from the context, or is irrelevant, we shall continue to simply write $T^*\U$. We point out that the zero-section of $T^*\U$ is the conormal space to $\U$ in $\U$, i.e., $T^*_{\U}\U$.

We will projectivize the fibers of $(T^*\U)^\C$ (resp., $(T^*\U)^\R$), and denote this projectivization by $\Proj\big((T^*\U)^\C\big)$ (resp., $\Proj\big((T^*\U)^\R\big)$), which is isomorphic to $\U\times\Proj^n$ (resp., $\U\times \R\Proj^{2n+1}$).  In either the complex or real case, we let $\hat\pi$ denote the projection from the projectivization of $T^*\U$ to $\U$, and if $\eta$ is a non-zero element of the fiber $(T^*\U)_p$, we denote its projective class by $[\eta]$. 

A subset $E\subseteq T^*\U$ is $\C$-conic (resp., $\R$-conic) if $(p, \eta)\in E$ implies that, for all $a\in\C$ (resp., $a\in\R$), $(p, a\eta)\in E$. If $E$ is any subset of $T^*\U$, we let $\Proj(E)$ denote the (real or complex) projectivization $\{(p, [\eta])\ |\ (p, \eta)\in E-T^*_{\U}\U\}$, and let $E_p:=\pi^{-1}(p)\cap E$.

\medskip

We need the following easy lemmas.

\begin{lem}\label{lem:closed} Suppose that $E\subseteq T^*\U$ is closed and $\R$-conic (resp., $\C$-conic). Then, $\Proj(E)$ is closed in $\Proj\big((T^*\U)^\R\big)$ (resp., $\Proj\big((T^*\U)^\C\big)$) and $\pi(E)$ is closed in $\U$.
\end{lem}

\begin{proof} We shall prove the real case. The proof over the complex numbers is the same. Throughout, we shall write simply $T^*\U$, in place of $(T^*\U)^\R$.

By definition of the quotient topology on $\Proj(T^*\U)$, $\Proj(E)$ is closed if and only if 
$$E^\prime:=\{(p, \eta)\in T^*\U-T^*_\U\U\ |\ (p, [\eta])\in\Proj(E)\}$$ is closed in 
$T^*\U-T^*_\U\U$. As $E$ is conic, $E^\prime= E-T^*_\U\U$, which is closed in $T^*\U-T^*_\U\U$, since $E$ is closed in $T^*\U$. Thus, $\Proj(E)$ is closed.

\smallskip

Now, suppose that we have a sequence $p_i\in\pi(E)$ and $p_i\rightarrow p\in\U$. We need to show that $p\in \pi(E)$. Let $\eta_i$ be such that $(p_i, \eta_i)\in E$. Identify $T^*\U$ with $\U\times \R^{2n+2}$.

If an infinite number of the $\eta_i$ are zero, then, by taking a subsequence (which we continue to write as $p_i$), we have an infinite sequence $(p_i, 0)\in E$. Then, $(p_i, 0)\rightarrow (p, 0)\in E$, as $E$ is closed. Thus, $p=\pi(p,0)\in\pi(E)$.

If an infinite number of the $\eta_i$ are not zero, we can take a subsequence $(p_i, \eta_i/|\eta_i|)$, which is still in $E$, as $E$ is conic. Since the $\eta_i/|\eta_i|$ are contained in the unit sphere, by taking another subsequence, we may assume that $\eta_i/|\eta_i|$ converges to some $\eta$. Thus, $(p_i, \eta_i/|\eta_i|)\rightarrow (p, \eta)$, which is in $E$, since $E$ is closed, and so $p\in\pi(E)$.
\end{proof}

\begin{lem}\label{lem:closedcrit} Suppose that $\tilde h_1, \dots, \tilde h_k$ are real (resp., complex) analytic functions from $\U$ to $\R$ (resp., $\C$), and suppose that $E\subseteq T^*\U$ is closed and $\R$-conic (resp., $\C$-conic). Then, the set ${\overline{\Sigma}}_E(\tilde h_1, \dots, \tilde h_k)$ of $p\in\U$ such that there exists  non-zero $(a_1, \dots, a_k)\in \R^k$ (resp., $\C^k$) such that 
$$
a_1d_p\tilde h_1+\dots+a_kd_p\tilde h_k\in E_p
$$
is closed in $\U$.
\end{lem}
\begin{proof} We shall prove the real case. The proof over the complex numbers is the same. Let $N:=2n+2$.

Let $K$ be the set of points $p\in\U$ such that $d_p\tilde h_1, \dots, d_p\tilde h_k$ are linearly dependent, i.e., let $K$ be the critical locus of the map $(\tilde h_1, \dots, \tilde h_k)$. Note that $K$ is closed.

Consider the continuous function $T:(\U-K)\times\R\Proj^{k-1}\rightarrow (\U-K)\times\R\Proj^{N-1}$ given by $T(p, [a_1,\dots, a_k]) = (p, [a_1d_p\tilde h_1+\dots+a_kd_p\tilde h_k])$. Let $\check\pi:(\U-K)\times\R\Proj^{k-1}\rightarrow\U-K$ denote the projection map.

By \lemref{lem:closed}, $\Proj(E)$ is closed in $\U\times\R\Proj^{N-1}$. Therefore, $B:=T^{-1}\big(((\U-K)\times\R\Proj^{N-1})\cap \Proj(E)\big)$ is closed in $(\U-K)\times\R\Proj^{k-1}$. As $\check\pi$ is proper, $\check\pi(B)$ is closed in $\U-K$.

Now, the set ${\overline{\Sigma}}_E(\tilde h_1, \dots, \tilde h_k)$ is equal to $\check\pi(B)\cup K$, which is closed in $\U$.
\end{proof}

\smallskip

\begin{defn} Let $\widetilde H:= (\tilde h_1, \dots, \tilde h_k)$ and $E$ be as in  \lemref{lem:closedcrit}. Then, the set ${\overline{\Sigma}}_E\widetilde H$ from \lemref{lem:closedcrit} is the {\bf closed $E$-critical locus of $\widetilde H$}. We define the {\bf $E$-discriminant of $\widetilde H$}, $\Delta_{E}\widetilde H$, to be $\widetilde H({\overline{\Sigma}}_{E}\widetilde H)$.

If $\widetilde H$ is complex analytic, and $E=SS(\Fdot)$, then we let ${\overline{\Sigma}}_{\Fdot}\widetilde H:= {\overline{\Sigma}}_E\widetilde H$ and $\Delta_{\Fdot}\widetilde H:=\Delta_{E}\widetilde H$.
\end{defn}

\smallskip

\begin{rem} By \propref{prop:fdotcrit}, if $E:=SS(\Fdot)$, then $\overline{\Sigma_{\Fdot}f}=\overline{\Sigma}_E\tilde f$; this was our reason for adopting our notation for the closed $E$-critical locus.

While we shall not need it in this paper, it is possible to show that, in special cases, there is a reasonable notion of the (non-closed) critical locus $\Sigma_E(h_1, \dots, h_k)$ which depends only on the restriction $(h_1, \dots, h_k)$ of $(\tilde h_1, \dots, \tilde h_k)$ to $\pi(E)$, and  $\overline{\Sigma_E(h_1, \dots, h_k)}= {\overline{\Sigma}}_E(\tilde h_1, \dots, \tilde h_k)$. In particular, this is the case when $\widetilde H$ is complex analytic and $E=SS(\Fdot)$.
\end{rem}

\smallskip

\begin{lem}\label{lem:lem} Suppose that $Z$ is an analytic subset of $\U$ such that $\tilde f$ is not constant on any irreducible component of $Z$. Let $M$ be an open dense subset of $Z_{\operatorname{reg}}$.  Then,
for all $p\in Z$, 
$$
\big(\overline{T^*_M\U}\big)_p + \langle d_p\tilde f\rangle:=\big\{\omega+ad_p\tilde f\ |\ a\in\C, \omega\in \big(\overline{T^*_M\U}\big)_p\big\} \subseteq \big(\overline{T^*_{\tilde f_{|_M}}\U}\big)_p,
$$
and, if $d_p\tilde f\not\in \big(\overline{T^*_M\U}\big)_p$, then this containment is an equality.
\end{lem}
\begin{proof} Since $\tilde f$ is not constant on any component of $Z$, $M^\prime:=M-\Sigma(\tilde f_{|_M})$ is dense in $Z$. Thus, $\overline{T^*_M\U}=\overline{T^*_{M^\prime}\U}$, and $\overline{T^*_{\tilde f_{|_M}}\U}=\overline{T^*_{\tilde f_{|_{M^\prime}}}\U}$.

Consider $\eta:=\omega+ad_p\tilde f\in \big(\overline{T^*_M\U}\big)_p + \langle d_p\tilde f\rangle$, where $(p_i, \omega_i)\in T^*_{M^\prime}\U$ and $(p_i, \omega_i)\rightarrow (p, \omega)$. Then, $\eta_i:=\omega_i+ad_{p_i}\tilde f\in \big(T^*_{\tilde f_{|_M}}\U\big)_{p_i}$, and $\eta_i\rightarrow\eta$. Therefore, the containment holds.

Suppose now that $d_p\tilde f\not\in \big(\overline{T^*_M\U}\big)_p$, and that $\eta\in \big(\overline{T^*_{\tilde f_{|_M}}\U}\big)_p=\big(\overline{T^*_{\tilde f_{|_{M^\prime}}}\U}\big)_p$. Then, there exists an analytic path $(p(t), \eta_t)\in \overline{T^*_{\tilde f_{|_{M^\prime}}}\U}$ such that $p=p(0)$, $\eta=\eta_0$, and, for $t\neq 0$, $(p(0), \eta_t)\in T^*_{\tilde f_{|_{M^\prime}}}\U$. Hence, for $t\neq 0$, $\eta_t=\omega_t+\lambda(t)d_{p(t)}\tilde f$, where $\lambda(t)\in\C$ and $\omega_t\in \big(T^*_M\U\big)_{p(t)}$, and $\lambda(t)$ and $\omega_t$ are uniquely determined. Evaluating at $p^\prime(t)$, we find that, for $t\neq 0$, $\eta_t(p^\prime(t))= \lambda(t)d_{p(t)}\tilde f(p^\prime(t))$, and so $\lambda(t)$ is a quotient of two analytic functions. Therefore, there are two possibilities as $t\rightarrow 0$: either $\lambda(t)$ approaches some $a\in\C$, or $|\lambda(t)|\rightarrow\infty$.

If $\lambda(t)\rightarrow a\in \C$, then $\omega_t\rightarrow \eta-ad_p\tilde f$, and $\eta=( \eta-ad_p\tilde f)+ ad_p\tilde f\in \big(\overline{T^*_M\U}\big)_p + \langle d_p\tilde f\rangle$. We claim that $d_p\tilde f\not\in \big(\overline{T^*_M\U}\big)_p$ implies that the case $|\lambda(t)|\rightarrow\infty$ cannot occur. Once we show this, the proof will be finished.

If $|\lambda(t)|\rightarrow\infty$, then, as $\eta_t\rightarrow\eta$,
$$
\frac{\omega_t}{\lambda(t)}+d_{p(t)}\tilde f = \frac{\eta_t}{\lambda(t)} \rightarrow 0,
$$
i.e., 
$$d_p\tilde f=\lim_{t\rightarrow 0}\left(-\frac{\omega_t}{\lambda(t)}\right)\in \big(\overline{T^*_M\U}\big)_p.$$
\end{proof}

\begin{lem} \label{lem:lemlem}
\begin{enumerate}
\item $SS_p(\Fdot)+\langle d_p\tilde f\rangle\subseteq \big|\big(T^*_{{}_{f,\Fdot}}\U\big)^\bullet\big|_p\cup SS_p(\Fdot)$.
\item Suppose that $p\not\in\overline{\Sigma_{\Fdot}f}$. Then, $\big|\big(T^*_{{}_{f,\Fdot}}\U\big)^\bullet\big|_p = SS_p(\Fdot)+\langle d_p\tilde f\rangle$.
\end{enumerate}
\end{lem}
\begin{proof} 

\noindent Proof of Item 1:

Suppose that $\eta\in SS_p(\Fdot)+\langle d_p\tilde f\rangle$. Then, $\eta =\omega+ad_p\tilde f$, where $a\in\C$ and $\omega\in\big(\overline{T^*_S\U}\big)_p$ for some $S\in\strat(\Fdot)$. Then, there exist $(p_i, \omega_i)\in T^*_S\U$ such that $(p_i, \omega_i)\rightarrow (p, \omega)$. There are two cases.

If $f$ is constant on $S$, then $(p_i, \omega_i+ad_{p_i}\tilde f)\in T^*_S\U$. Hence, $\eta=\omega+ad_p\tilde f\in \big(\overline{T^*_S\U}\big)_p\subseteq SS_p(\Fdot)$. If $f$ is not constant on $S$, then \lemref{lem:lem} implies that $\eta\in \big|\big(T^*_{{}_{f,\Fdot}}\U\big)^\bullet\big|_p$.

\medskip

\noindent Proof of Item 2:

Note that if $S\in\strat(\Fdot)$ and $f_{|_S}$ is constant, then \propref{prop:fdotcrit} implies that $S\subseteq\overline{\Sigma_{\Fdot}f}$; hence, by our hypothesis, $p\not\in\overline{S}$. Therefore,
$$
SS_p(\Fdot)= \pi^{-1}(p)\cap\bigcup_{\substack{S\in\strat(\Fdot)\\ f_{|_S}\neq{\rm\ const.}}}\overline{T^*_S\U}.
$$
Now the result follows immediately from \lemref{lem:lem} and the definition of $\big(T^*_{{}_{f,\Fdot}}\U\big)^\bullet$.
\end{proof}

\smallskip

\begin{prop}\label{prop:sigma} There is an equality of sets
$${\overline{\Sigma}}_{\Fdot}(\tilde f, \tilde g)\ =\ \overline{\Sigma_{\Fdot}f}\ \cup\ \big|\Gamma_{f, \tilde g}(\Fdot)\big|.$$
\end{prop}
\begin{proof} Suppose that $p\not\in\overline{\Sigma_{\Fdot}f}$. Then, it follows immediately from Item 2 of \lemref{lem:lemlem} that $p\in {\overline{\Sigma}}_{\Fdot}(\tilde f, \tilde g)$ if and only if $p\in \big|\Gamma_{f, \tilde g}(\Fdot)\big|$. Therefore,
$$
{\overline{\Sigma}}_{\Fdot}(\tilde f, \tilde g)- \overline{\Sigma_{\Fdot}f}\ =\ \big|\Gamma_{f, \tilde g}(\Fdot)\big|-\overline{\Sigma_{\Fdot}f}.
$$
Now, take the union of both sides above with $\overline{\Sigma_{\Fdot}f}$, and use that $\overline{\Sigma_{\Fdot}f}\subseteq {\overline{\Sigma}}_{\Fdot}(\tilde f, \tilde g)$.
\end{proof}

\smallskip

\begin{rem}\label{rem:symm2}  Note that is no claim in \propref{prop:sigma} about $\overline{\Sigma_{\Fdot}f}$ and  $\big|\Gamma_{f, \tilde g}(\Fdot)\big|$ intersecting in some nice way. In fact, in \remref{rem:main2}, we give an example where these two sets are equal. However, if   $\big|\Gamma_{f, \tilde g}(\Fdot)\big|$ is $1$-dimensional and $f$ is not constant along any component of $\big|\Gamma_{f, \tilde g}(\Fdot)\big|$, then \propref{prop:fdotcrit} implies that the intersection of $\overline{\Sigma_{\Fdot}f}$ and  $\big|\Gamma_{f, \tilde g}(\Fdot)\big|$ is either empty or consists of isolated points.

We also want to return to the topic of symmetry that we first discussed in \remref{rem:polar}. By \propref{prop:sigma} and the symmetry of the definition of the closed critical locus, we have that 
$$
\overline{\Sigma_{\Fdot}f}\ \cup\ \big|\Gamma_{f, \tilde g}(\Fdot)\big| = \overline{\Sigma_{\Fdot}g}\ \cup\ \big|\Gamma_{g, \tilde f}(\Fdot)\big|.
$$
By \propref{prop:fdotcrit}, $f$ and $g$ are constant along the components of $\overline{\Sigma_{\Fdot}f}$ and $\overline{\Sigma_{\Fdot}g}$, respectively. It follows that the symmetric components of $\big|\Gamma_{f, \tilde g}(\Fdot)\big|$ and $\big|\Gamma_{g, \tilde f}(\Fdot)\big|$ are the same.

Note, however, that even in the classical case where we look at germs at the origin, $f$ is fixed, and $\tilde g$ is chosen to be a generic linear form, it is, in general, {\bf false} that there is an equality of sets $\big|\Gamma_{f, \tilde g}(\Fdot)\big|=\big|\Gamma_{g, \tilde f}(\Fdot)\big|$. Suppose, for instance, that  $\dim_\0\overline{\Sigma_{\Fdot}f}\geq 2$. For a generic linear form $\mathfrak l$, either $\overline{\Sigma_{\Fdot}\mathfrak l}$ will be empty or the origin will be an isolated point in $\overline{\Sigma_{\Fdot}\mathfrak l}$; furthermore, \propref{prop:geng} implies that  $\big|\Gamma_{f, \tilde g}(\Fdot)\big|$ is purely $1$-dimensional at the origin. However,  
$$\big|\Gamma_{g, \tilde f}(\Fdot)\big| \  = \  \overline{\Sigma_{\Fdot}g}\ \cup\ \big|\Gamma_{g, \tilde f}(\Fdot)\big| \  =  \ \overline{\Sigma_{\Fdot}f}\ \cup\ \big|\Gamma_{f, \tilde g}(\Fdot)\big|
$$
is, at least, $2$-dimensional at the origin.
\end{rem}

\bigskip

We now need to prove our main technical lemma.

\smallskip

Let ${\stackrel{\circ}{\D}}_{{}_\delta}$ denote an open disk of radius $\delta$, centered at the origin, in $\C$. For positive $\rho, \delta, \epsilon\in\R$, let $N^\epsilon_{\delta, \rho}:=B_\epsilon\cap g^{-1}({\stackrel{\circ}{\D}}_{{}_\delta})\cap f^{-1}({\stackrel{\circ}{\D}}_{{}_\rho})$, let $P^\epsilon_{\delta, \rho}$ be the restriction of $(\tilde g, \tilde f)$ to a map from $N^\epsilon_{\delta, \rho}$ to ${\stackrel{\circ}{\D}}_{{}_\delta}\times {\stackrel{\circ}{\D}}_{{}_\rho}$, and let $(\Fdot)^\epsilon_{\delta, \rho}$ denote the restriction of $\Fdot$ to $N^\epsilon_{\delta, \rho}$. Let $r:\U\rightarrow\R$ be the ``squared distance from the origin'' function $r=|z_0|^2+\dots+|z_n|^2$.

\begin{lem}\label{lem:main2}  Suppose that $\dim_\0 V(f)\cap \big|\Gamma_{f, \tilde g}(\Fdot)\big|\leq 0$. Then, there exists $\epsilon_0>0$ such that, for all $\epsilon_1$ and $\epsilon_2$ such that $0<\epsilon_2<\epsilon_1\leq\epsilon_0$, there exist $\delta, \rho>0$ such that
\begin{enumerate}
\item for all $p\in N^{\epsilon_1}_{\delta, \rho}-{\stackrel{\circ}{B}}_{\epsilon_2}$, 
$$
d_pr\not\in \big(SS_p(\Fdot)+<d_p\tilde g, d_p\tilde f>\big)^\R;
$$
\item for all $p\in \big(N^{\epsilon_1}_{\delta, \rho}-{\stackrel{\circ}{B}}_{\epsilon_2}\big)-\overline{\Sigma_{\Fdot}f}$, in particular,  for all $p\in \big(N^{\epsilon_1}_{\delta, \rho}-{\stackrel{\circ}{B}}_{\epsilon_2}\big)-V(f)$,  for all $(a,b)\in\C^2-\{\0\}$,
$$\operatorname{Re}(ad_p\tilde f+bd_p\tilde g)\not\in (SS_p(\Fdot))^\R+\langle d_pr\rangle^\R.$$
\end{enumerate}
\end{lem}
\begin{proof}

\noindent Proof of Item 1:

\smallskip

Let $Y$ be the set of $p\in U$ such that there exists non-zero $(c, b)\in\R\times\C$ such that $cd_pr+\operatorname{Re}(bd_p\tilde g)\in \big|\big(T^*_{{}_{f,\Fdot}}\U\big)^\bullet\big|_p^\R$.  By \lemref{lem:closedcrit}, $Y$ is closed in $\U$. Let $Z$ be the set of $p\in U$ such that 
$$d_pr\in \big(SS_p(\Fdot)+<d_p\tilde g, d_p\tilde f>\big)^\R.$$
Let $B_{\epsilon_0}^*:=B_{\epsilon_0}-\{\0\}$.

We shall prove Item 1 by proving that there exists $\epsilon_0>0$ such that 
\begin{enumerate}\item[a.] $B_{\epsilon_0}^*\cap Y\cap V(f, g) = \emptyset$;
\item[b.] $B_{\epsilon_0}^*\cap Z- V(f) \subseteq B_{\epsilon_0}^*\cap Y-V(f)$; and
\item[c.] for all $\epsilon_1, \epsilon_2$ such that $0<\epsilon_2<\epsilon_1\leq\epsilon_0$, there exists $\delta^\prime>0$ such that 
$$(B_{\epsilon_1}-{\stackrel{\circ}{B}}_{\epsilon_2})\cap Z\cap g^{-1}({\stackrel{\circ}{\D}}_{{}_{\delta^\prime}})\cap V(f)=\emptyset.$$
\end{enumerate}

We will first show how Items a, b, and c imply Item 1. We will then show that Items a, b, and c hold.

Assume Items a, b, and c, and let $\epsilon_1, \epsilon_2$ be such that $0<\epsilon_2<\epsilon_1\leq\epsilon_0$. Let $\delta^\prime$ be as in Item c. As $Y$ is closed, $(B_{\epsilon_1}-{\stackrel{\circ}{B}}_{\epsilon_2})\cap Y$ is compact. Thus, since $(B_{\epsilon_1}-{\stackrel{\circ}{B}}_{\epsilon_2})\cap Y\cap V(f, g)=\emptyset$ by Item a, there exist $\rho, \delta^{\prime\prime}>0$ such that  
$$(B_{\epsilon_1}-{\stackrel{\circ}{B}}_{\epsilon_2})\cap Y \cap g^{-1}({\stackrel{\circ}{\D}}_{{}_{\delta^{\prime\prime}}})\cap f^{-1}({\stackrel{\circ}{\D}}_{{}_\rho})=\emptyset.\leqno(\dagger)$$
Fix such  $\rho$ and $\delta^{\prime\prime}$. Fix $\delta$ such that  $0<\delta\leq\operatorname{min}\{\delta^{\prime}, \delta^{\prime\prime}\}$. We wish to show that $\big(N^{\epsilon_1}_{\delta, \rho}-{\stackrel{\circ}{B}}_{\epsilon_2}\big)\cap Z=\emptyset$. Suppose that $p\in \big(N^{\epsilon_1}_{\delta, \rho}-{\stackrel{\circ}{B}}_{\epsilon_2}\big)\cap Z$. Then, $p\not\in V(f)$ by Item c. However, then, Item b implies that $p\in Y$; a contradiction of $(\dagger)$.

\bigskip

Now we will show that we may pick  $\epsilon_0>0$ so that Items a, b, and c hold. 
Choose $\epsilon_0>0$ such that, for all $\epsilon$ such that $0<\epsilon\leq\epsilon_0$, $\partial B_\epsilon$ transversely intersects all of the strata of $\strat$. Then, for all $p\in B_{\epsilon_0}-\{\0\}$, $d_pr\not\in (SS_p(\Fdot))^\R$. Similarly, we may also choose $\epsilon_0>0$ so that, for all $\epsilon$ such that $0<\epsilon\leq\epsilon_0$, for all $p\in (B_{\epsilon_0}-\{\0\})$, 
$$d_pr\not\in (SS_p(\psi_f[-1]\Fdot))^\R\cup (SS_p(\psi_g[-1]\Fdot))^\R\cup (SS_p(\psi_g[-1]\psi_f[-1]\Fdot))^\R.$$
Combining the paragraph above with \propref{prop:fdotcrit} and the equivalences at the beginning of \thmref{thm:main1}, and using our hypothesis that $\dim_\0 V(f)\cap \big|\Gamma_{f, \tilde g}(\Fdot)\big|\leq 0$, we may pick $\epsilon_0>0$ such that
\begin{enumerate}
\item[i.] $B^*_{\epsilon_0}\cap {\overline{\Sigma}}_{\Fdot}f\subseteq V(f)$;
\item[ii.] $B^*_{\epsilon_0}\cap {\overline{\Sigma}}_{{}_{\psi_f[-1]\Fdot}}g=\emptyset$;
\item[iii.] for all $p\in B^*_{\epsilon_0}$, 
$$d_pr\not\in (SS_p(\Fdot))^\R\cup(SS_p(\psi_f[-1]\Fdot))^\R\cup (SS_p(\psi_g[-1]\Fdot))^\R\cup (SS_p(\psi_g[-1]\psi_f[-1]\Fdot))^\R.$$
\end{enumerate}

\bigskip

\noindent Proof of Item a:

Suppose that we have non-zero $(c, b)\in\R\times \C$ and $p\in B^*_{\epsilon_0}\cap V(f,g)$ such that $cd_pr+\operatorname{Re}(bd_p\tilde g)\in \big|\big(T^*_{{}_{f,\Fdot}}\U\big)^\bullet\big|^\R_p=(SS_p(\psi_f[-1]\Fdot))^\R$, where the last equality follows from \thmref{thm:psigecc}. If $c=0$, then $d_p\tilde g\in SS_p(\psi_f[-1]\Fdot)$, i.e., $p\in {\overline{\Sigma}}_{{}_{\psi_f[-1]\Fdot}}g$, which contradicts Item ii. Thus, $c$ must be unequal to zero, and so $d_pr\in (SS_p(\psi_f[-1]\Fdot)+<d_p\tilde g>)^\R$. As $d_p\tilde g\not\in SS_p(\psi_f[-1]\Fdot)$ and $p\in V(g)$, by \lemref{lem:lemlem}, Item 2, $SS_p(\psi_f[-1]\Fdot)+<d_p\tilde g>= SS_p(\psi_g[-1]\psi_f[-1]\Fdot)$. However, $d_pr\in SS_p(\psi_g[-1]\psi_f[-1]\Fdot)$ contradicts Item iii. This proves Item a.

\bigskip

\noindent Proof of Item b:

Suppose that $p\in B_{\epsilon_0}^*\cap Z- V(f)$, i.e., $p\in B_{\epsilon_0}^*- V(f)$ and there exists $b\in \C$ such that $d_pr +\operatorname{Re}(bd_p\tilde g)\in \big(SS_p(\Fdot)+<d_p\tilde f>\big)^\R$. By Item i and \lemref{lem:lemlem}, Item 2, $SS_p(\Fdot)+<d_p\tilde f> =  \big|\big(T^*_{{}_{f,\Fdot}}\U\big)^\bullet\big|_p$, and so $d_pr +\operatorname{Re}(bd_p\tilde g)\in  \big|\big(T^*_{{}_{f,\Fdot}}\U\big)^\bullet\big|_p^\R$. Hence, $p\in Y$, and we have proved Item b.

\bigskip

\noindent Proof of Item c:

Let $W$ be the closed set 
$$
W:=\{p\in\U\ |\ d_pr\in\big( \big|\big(T^*_{{}_{g,\psi_f[-1]\Fdot\Fdot}}\U\big)^\bullet\big|_p\cup  \big|\big(T^*_{{}_{g,\Fdot}}\U\big)^\bullet\big|_p\cup SS_p(\psi_f[-1]\Fdot)\cup SS_p(\Fdot)\big)^\R\}.
$$
By Item iii and \thmref{thm:psigecc}, $B^*_{\epsilon_0}\cap W\cap V(g)=\emptyset$. We also claim that
$$
B^*_{\epsilon_0}\cap Z\cap V(f)\subseteq B^*_{\epsilon_0}\cap W\cap V(f).\leqno(\ddagger)
$$
For $p\in Z$ if and only if there exists a non-zero $b\in\C$ such that $d_pr+\operatorname{Re}(bd_p\tilde g)\in (SS_p(\Fdot)+<d_p\tilde f>)^\R$, and, by Item 1 of \lemref{lem:lemlem},
$$
SS_p(\Fdot)+<d_p\tilde f>\subseteq \big|\big(T^*_{{}_{f,\Fdot}}\U\big)^\bullet\big|_p\cup SS_p(\Fdot).
$$
If $p\in V(f)$, then \thmref{thm:psigecc} implies that $\big|\big(T^*_{{}_{f,\Fdot}}\U\big)^\bullet\big|_p = SS_p(\psi_f[-1]\Fdot)$. Now, $(\ddagger)$ follows at once.

Let $\epsilon_1, \epsilon_2$ be such that $0<\epsilon_2<\epsilon_1\leq\epsilon_0$. Then, since  $B^*_{\epsilon_0}\cap W\cap V(g)=\emptyset$ and $W$ is closed, $(B_{\epsilon_1}-{\stackrel{\circ}{B}}_{\epsilon_2})\cap W$ is compact and $(B_{\epsilon_1}-{\stackrel{\circ}{B}}_{\epsilon_2})\cap W\cap V(g)=\emptyset$. Thus, there exists $\delta^\prime>0$ such that $(B_{\epsilon_1}-{\stackrel{\circ}{B}}_{\epsilon_2})\cap W\cap g^{-1}({\stackrel{\circ}{\D}}_{{}_{\delta^\prime}})=\emptyset$. Finally, we conclude that
$$
(B_{\epsilon_1}-{\stackrel{\circ}{B}}_{\epsilon_2})\cap Z\cap g^{-1}({\stackrel{\circ}{\D}}_{{}_{\delta^\prime}})\cap V(f)\ \subseteq\ (B_{\epsilon_1}-{\stackrel{\circ}{B}}_{\epsilon_2})\cap W\cap g^{-1}({\stackrel{\circ}{\D}}_{{}_{\delta^\prime}})\cap V(f) = \emptyset,
$$
which proves Item c, and concludes the proof of Item 1 from the statement of the lemma.

\bigskip

\noindent Proof of Item 2:

\smallskip

Assume that $\dim_\0 V(f)\cap \big|\Gamma_{f, \tilde g}(\Fdot)\big|\leq 0$. Pick $\epsilon_0>0$ so that $B_{\epsilon_0}\cap V(f)\cap \big|\Gamma_{f, \tilde g}(\Fdot)\big|\subseteq\{\0\}$. This means precisely that 
$$
(B_{\epsilon_0}-\{\0\})\cap\{p\in \U\ |\ d_p\tilde g\in \big|\big(T^*_{{}_{f,\Fdot}}\U\big)^\bullet\big|_p\}\cap V(f)=\emptyset.\leqno(*)
$$
Let $\epsilon_1$ and $\epsilon_2$ be such that $0<\epsilon_2<\epsilon_1\leq\epsilon_0$, and assume that $\rho$ and $\delta$ are such that Item 1 holds. Then, 
$$(B_{\epsilon_1}-{\stackrel{\circ}{B}}_{\epsilon_2})\cap V(f)\cap \big|\Gamma_{f, \tilde g}(\Fdot)\big| =\emptyset.
$$
As $(B_{\epsilon_1}-{\stackrel{\circ}{B}}_{\epsilon_2})\cap \big|\Gamma_{f, \tilde g}(\Fdot)\big|$ is compact, it follows that we may re-choose $\rho$, smaller if needed, so that 
$$(B_{\epsilon_1}-{\stackrel{\circ}{B}}_{\epsilon_2})\cap f^{-1}({\stackrel{\circ}{\D}}_{{}_\rho})\cap \big|\Gamma_{f, \tilde g}(\Fdot)\big| =\emptyset. \leqno(\dagger)
$$
We claim that Item 2 holds.

To see this, let $p\in \big(N^{\epsilon_1}_{\delta, \rho}-{\stackrel{\circ}{B}}_{\epsilon_2}\big)-\overline{\Sigma_{\Fdot}f}$, and assume that we have $(a,b)\in\C^2-\{\0\}$ such that
$$\operatorname{Re}(ad_p\tilde f+bd_p\tilde g)\in (SS_p(\Fdot))^\R+\langle d_pr\rangle^\R.$$
Then, there exists a real number $c$ such that 
$$
cd_pr+ \operatorname{Re}(ad_p\tilde f+bd_p\tilde g)\in (SS_p(\Fdot))^\R.
$$
If $c\neq 0$, we may divide by $c$ and obtain a contradiction to Item 1. Thus, $c$ must equal $0$, and so 
$$
ad_p\tilde f+bd_p\tilde g\in SS_p(\Fdot),
$$
i.e., $p\in {\overline{\Sigma}}_{\Fdot}(\tilde f, \tilde g)$, which by \propref{prop:sigma}, is equal to $\overline{\Sigma_{\Fdot}f}\ \cup\ \big|\Gamma_{f, \tilde g}(\Fdot)\big|$. This is a contradiction of $(\dagger)$ and the fact that $p\not\in\overline{\Sigma_{\Fdot}f}$. \end{proof}

\bigskip

As before, for $\epsilon, \delta, \rho>0$, let $N^\epsilon_{\delta, \rho}:=B_\epsilon\cap g^{-1}({\stackrel{\circ}{\D}}_{{}_\delta})\cap f^{-1}({\stackrel{\circ}{\D}}_{{}_\rho})$. Let $(\Fdot)^\epsilon_{\delta, \rho}$ be the restriction of $\Fdot$ to $N^\epsilon_{\delta, \rho}$, and let $T^\epsilon_{\delta, \rho}$ be the restriction of the map $(g, f)$ to a map from $N^\epsilon_{\delta, \rho}$ to ${\stackrel{\circ}{\D}}_{{}_\delta}\times {\stackrel{\circ}{\D}}_{{}_\rho}$.

\smallskip

\begin{thm}\label{thm:main2} {\rm (The Derived Category Discriminant Theorem)} Suppose that $\dim_\0 V(f)\cap \big|\Gamma_{f, \tilde g}(\Fdot)\big|\leq 0$. 

Then, for all sufficiently small $\epsilon>0$, there exist $\delta, \rho>0$ such that the derived push-forward $\Adot:=R(T^\epsilon_{\delta, \rho})_*(\Fdot)^\epsilon_{\delta, \rho}$ is complex analytically constructible with respect to the stratification given by 
$$
\{{\stackrel{\circ}{\D}}_{{}_\delta}\times {\stackrel{\circ}{\D}}_{{}_\rho}-\Delta_{\Fdot}(\tilde g, \tilde f),\ ({\stackrel{\circ}{\D}}_{{}_\delta}\times {\stackrel{\circ}{\D}}_{{}_\rho})\cap\Delta_{\Fdot}(\tilde g, \tilde f)-\{\0\}, \ \{\0\}\}.
$$
\end{thm}
\begin{proof} Fix choices of $\epsilon_0$, $\epsilon_1$, $\epsilon_2$, $\delta$, and $\rho$ as in \lemref{lem:main2}.  Pick $\epsilon$ so that $\epsilon_2<\epsilon<\epsilon_1$. Let  $\Gdot$ be the restriction of $\Fdot$ to $Y:=\stackrel{\circ}{B}_{\epsilon_1}\cap g^{-1}({\stackrel{\circ}{\D}}_{{}_\delta})\cap f^{-1}({\stackrel{\circ}{\D}}_{{}_\rho})$. Let $\tilde h$ be the restriction of $(g, f)$ to a map from $Y$ to ${\stackrel{\circ}{\D}}_{{}_\delta}\times {\stackrel{\circ}{\D}}_{{}_\rho}$. Then, by applying Item 1 of \lemref{lem:main2} and Proposition 8.5.8 of \cite{kashsch} to $\Gdot$ (where the $\phi$ and $f$ of \cite{kashsch} are our $r$ and $(g, f)$, respectively),  we immediately conclude that $\Adot$ is complex analytically constructible. 

As $\Delta_{\Fdot}(\tilde g, \tilde f)$ is either empty or a curve, to show that $\Adot$ is constructible with respect to the given stratification, one has only to show that the cohomology of $\Adot$ is locally constant at points in ${\stackrel{\circ}{\D}}_{{}_\delta}\times {\stackrel{\circ}{\D}}_{{}_\rho}-\Delta_{\Fdot}(\tilde g, \tilde f)$. Let  $q:=(u_0, v_0)\in {\stackrel{\circ}{\D}}_{{}_\delta}\times {\stackrel{\circ}{\D}}_{{}_\rho}$; it suffices to show that $SS_q(\Adot)=\{0\}$.

In the following, we use the real structure in each of the statements. Proposition 8.5.8 of \cite{kashsch} implies Proposition 5.4.17 of \cite{kashsch}. Item ii, part d, of this latter proposition tells us that 
$$SS_q(\Adot)\subseteq\bigcup_{p\in B_\epsilon\cap g^{-1}(u_0)\cap f^{-1}(v_0)}\{ad_p\tilde g+bd_p\tilde f\in SS_p((\Fdot)^\epsilon_{\delta, \rho})\ |\ (a,b)\in\C^2\}.
$$
If $p\in\stackrel{\circ}{B}_\epsilon$, then $SS_p((\Fdot)^\epsilon_{\delta, \rho}) = SS_p(\Fdot)$, and so $\{ad_p\tilde g+bd_p\tilde f\in SS_p((\Fdot)^\epsilon_{\delta, \rho})\ |\ (a,b)\in\C^2\}=\{0\}$ as $p\not\in\overline{\Sigma}_{\Fdot}(\tilde f, \tilde g)$. If $p\in\partial B_\epsilon$, then, by Proposition 5.4.8 of \cite{kashsch}, $SS_p((\Fdot)^\epsilon_{\delta, \rho}) \subseteq SS_p(\Fdot)+<d_pr>$, and so $\{ad_p\tilde g+bd_p\tilde f\in SS_p((\Fdot)^\epsilon_{\delta, \rho})\ |\ (a,b)\in\C^2\}=\{0\}$ by Item 2 of \lemref{lem:main2}. 
\end{proof}

\medskip

We refer to $\Adot$ in \thmref{thm:main2}, for $\epsilon$, $\delta$, and $\rho$ as in the theorem, as the {\it derived category discriminant of $(g, f)$}\/ or as the  {\it discriminant of $(g, f)$ as a complex of sheaves}.

\medskip

\begin{rem}\label{rem:main2} The assumption that $\dim_\0 V(f)\cap \big|\Gamma_{f, \tilde g}(\Fdot)\big|\leq 0$ is crucial in \lemref{lem:main2} and \thmref{thm:main2}. Consider the classic example of the map $H:=(\tilde g, \tilde f) = (g, f):\C^3\rightarrow\C^2$ given by $g(x,y,t)=x$ and $f(x,y,t)=y^2-tx^2$, where it is not possible to stratify the domain and codomain in order to obtain a Thom map. The (ordinary) discriminant of $H$ is simply the origin and, yet, for $0<\delta, \rho\ll\epsilon\ll1$, the isomorphism-type of the cohomology of the fibers $B_\epsilon\cap H^{-1}(a,b)$ is not independent of the choice of $(a,b)\in {\stackrel{\circ}{\D}}_{{}_\delta}\times {\stackrel{\circ}{\D}}_{{}_\rho}-\{\0\}$.

The reader should verify that, in this example, $\big|\Gamma_{f, \tilde g}(\Fdot)\big| = V(x, y)$ and so the condition that $\dim_\0 V(f)\cap \big|\Gamma_{f, \tilde g}(\Fdot)\big|\leq 0$  does not hold.
\end{rem}

\smallskip

\begin{cor}\label{cor:main2} Suppose that $\dim_\0 V(f)\cap \big|\Gamma_{f, \tilde g}(\Fdot)\big|\leq 0$. Let $\big(\widehat\Gamma^1_{f, \tilde g}(\Fdot)\big)^\bullet$ denote the components of $\big(\Gamma^1_{f, \tilde g}(\Fdot)\big)^\bullet$ which are not contained in $V(g)$. Then,
\begin{enumerate}
\item
$$
\hyp^{k-1}(F_{f, \mathbf 0}, F_{f_{|_{V(g)}}, \mathbf 0} ; \Fdot) \cong \big(\big(\widehat\Gamma^1_{f, \tilde g}(\Fdot)\big)^k\odot V(f)\big)_\0;
$$

\item
$$
\hyp^{k-1}(F_{g, \mathbf 0}, F_{g_{|_{V(f)}}, \mathbf 0} ; \Fdot) \cong \big(\big(\widehat\Gamma^1_{f, \tilde g}(\Fdot)\big)^k\odot V(g)\big)_\0; and
$$

\item
$$
H^k(\phi_f[-1]\psi_g[-1]\Fdot)_\0\ \cong\ \big(\big(\widehat\Gamma^1_{f, \tilde g}(\Fdot)\big)^k\odot V(g)\big)_\0\ \oplus\ H^k(\psi_g[-1]\phi_f[-1]\Fdot)_\0.
$$
\end{enumerate}
\end{cor}
\begin{proof} Now that we have \thmref{thm:main2}, the proof of each item is obtained by looking at the relative hypercohomology of a complex disk modulo a point, and using that this relative hypercohomology splits as a direct sum. One ``sees'' the results by looking at ``pictures'' in ${\stackrel{\circ}{\D}}_{{}_\delta}\times {\stackrel{\circ}{\D}}_{{}_\rho}$; exactly as in the case where $\Fdot$ is the constant sheaf on affine space and $g$ is a generic linear form. The discriminant/Cerf diagram arguments remain the same, {\bf except} that it is no longer true that the components of the Cerf diagram are tangent to the horizontal axis at the origin, i.e., it is not necessarily true for each component $C$ of $\big|\big(\widehat\Gamma_{f, \tilde g}(\Fdot)\big)^\bullet\big|$ that $(C\cdot V(f))_\0> (C\cdot V(g))_\0$.

Of course, the pictures are actually drawn in $\R^2$, and so a line segment represents a complex disk (but a point still represents a point). The three relevant pictures, in order, are:

\begin{center}{\includegraphics[height=30mm, width =60mm]{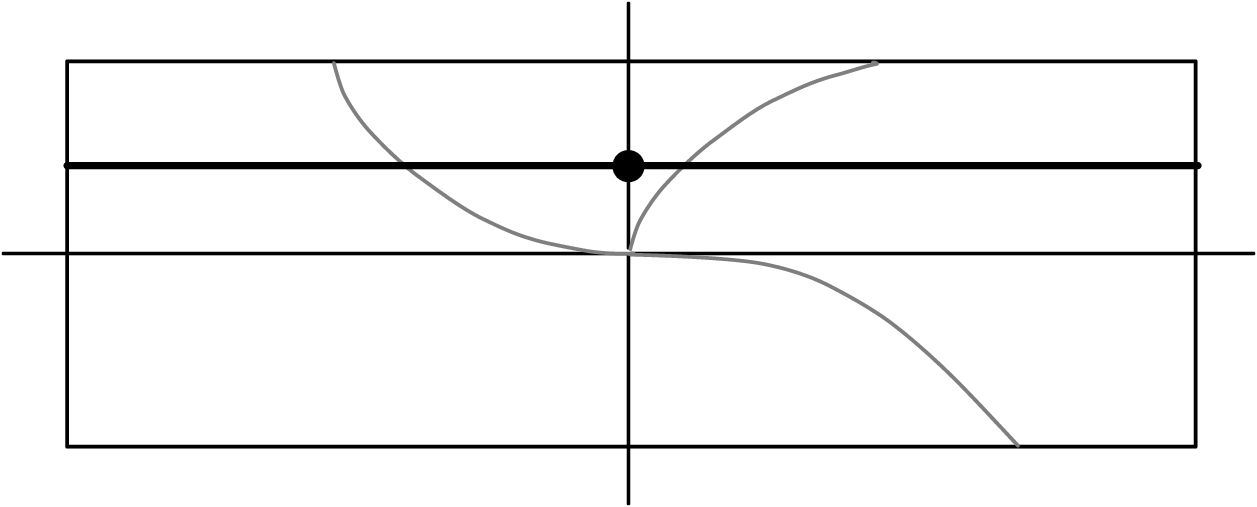}}\end{center}

\begin{center}{\includegraphics[height=30mm, width =60mm]{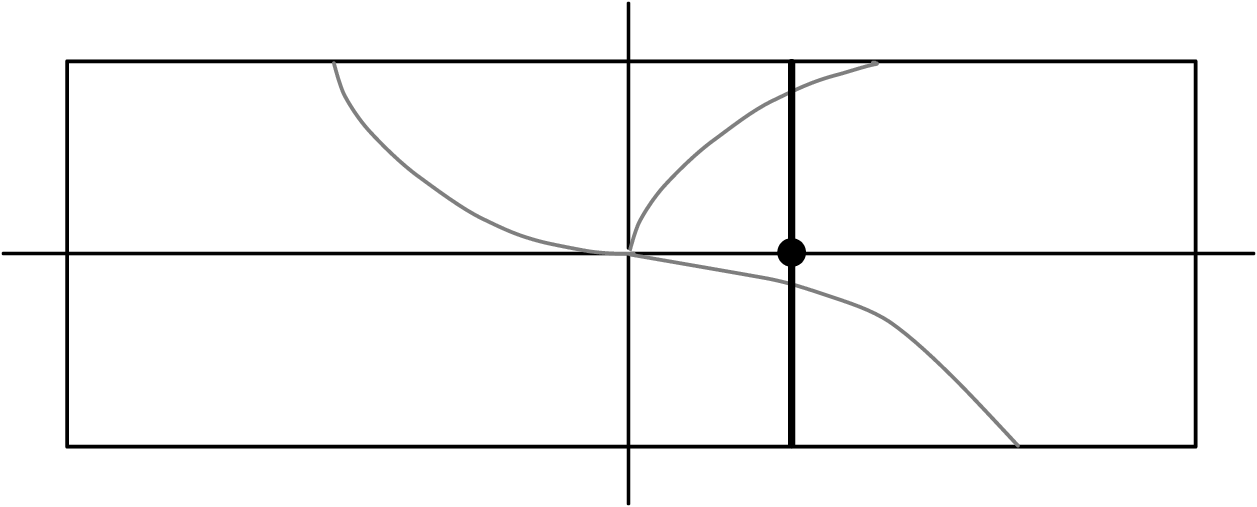}}\end{center}

\begin{center}{\includegraphics[height=30mm, width =60mm]{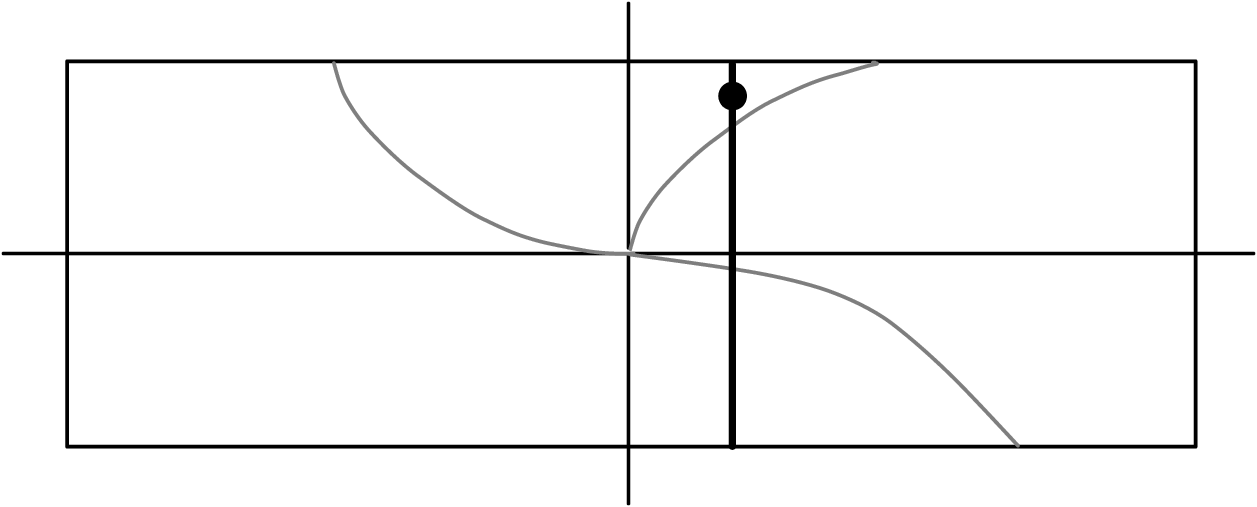}}\end{center}

\end{proof}

\bigskip

\begin{rem} If $\dim_\0 V(g)\cap \big|\Gamma_{f, \tilde g}(\Fdot)\big|\leq 0$, then $\big(\widehat\Gamma_{f, \tilde g}(\Fdot)\big)^\bullet=\big(\Gamma_{f, \tilde g}(\Fdot)\big)^\bullet$, and Item 1 of \corref{cor:main2}, combined with \thmref{thm:main1}, yields an isomorphism between the  cohomology $\hyp^{*}(F_{f, \mathbf 0}, F_{f_{|_{V(g)}}, \mathbf 0} ; \Fdot[-1])$ and $H^*(\phi_g[-1]\psi_f[-1]\Fdot)_\0$; this isomorphism is {\bf not} natural. In particular, the Milnor monodromy of $f$ typically induces completely different automorphisms on these two cohomologies.
\end{rem} 

\begin{cor}\label{cor:emptyequiv} The following are equivalent:
\begin{enumerate}
\item $\0\not\in\big|\Gamma_{f, \tilde g}(\Fdot)\big|$;
\item for all $p\in V(f, g)$ near the origin, $H^*(\phi_g[-1]\psi_f[-1]\Fdot)_p=0$;
\item $\dim_\0 V(f)\cap\big|\Gamma_{f, \tilde g}(\Fdot)\big|\leq 0$, and $H^*(\phi_g[-1]\psi_f[-1]\Fdot)_\0=0$;
\item $\dim_\0 V(g)\cap \big|\Gamma_{f, \tilde g}(\Fdot)\big|\leq 0$, and inclusion induces an isomorphism between $\hyp^{*}(F_{f, \mathbf 0} ; \Fdot)$ and $\hyp^{*}(F_{f_{|_{V(g)}}, \mathbf 0} ; \Fdot)$;
\item $\dim_\0 V(g)\cap \big|\Gamma_{f, \tilde g}(\Fdot)\big|\leq 0$, and inclusion induces an isomorphism between $\hyp^{*}(F_{g, \mathbf 0} ; \Fdot)$ and $\hyp^{*}(F_{g_{|_{V(f)}}, \mathbf 0} ; \Fdot)$;

\item $\dim_\0 V(g)\cap \big|\Gamma_{f, \tilde g}(\Fdot)\big|\leq 0$, and $
H^*(\phi_f[-1]\psi_g[-1]\Fdot)_\0$ is isomorphic to $H^*(\psi_g[-1]\phi_f[-1]\Fdot)_\0$.
\end{enumerate}
\end{cor}
\begin{proof} This is immediate from \lemref{lem:phipsi}, \thmref{thm:main1}, and \corref{cor:main2}.
\end{proof}

\smallskip

The above corollary has a more familiar feel in the classical case, where $\tilde g$ is a generic linear form. Recall that, in \propref{prop:geng}, we showed that, if  $\mathfrak l$ is a generic linear form on $\U$, then, $\dim_\0 V(f)\cap \big|\Gamma_{f, \mathfrak l}(\Fdot)\big|\leq 0$ and $\dim_\0 V(\mathfrak l)\cap \big|\Gamma_{f, \mathfrak l}(\Fdot)\big|\leq 0$. Thus, \corref{cor:emptyequiv} immediately yields:

\begin{cor}\label{cor:emptyequiv2} The following are equivalent:
\begin{enumerate}
\item there exists a non-zero linear form $\mathfrak l$ such that $\0\not\in\big|\Gamma_{f, \mathfrak l}(\Fdot)\big|$;

\item for generic linear $\mathfrak l$, $\0\not\in\big|\Gamma_{f, \mathfrak l}(\Fdot)\big|$;

\item for generic linear $\mathfrak l$, $H^*(\phi_{\mathfrak l}[-1]\psi_f[-1]\Fdot)_\0=0$;

\item for generic linear $\mathfrak l$, inclusion induces an isomorphism between $\hyp^{*}(F_{f, \mathbf 0} ; \Fdot)$ and $\hyp^{*}(F_{f_{|_{V(\mathfrak l)}}, \mathbf 0} ; \Fdot)$;

\item for generic linear $\mathfrak l$, inclusion induces an isomorphism between $\hyp^{*}(F_{\mathfrak l, \mathbf 0} ; \Fdot)$ and $\hyp^{*}(F_{{\mathfrak l}_{|_{V(f)}}, \mathbf 0} ; \Fdot)$;

\item for generic linear $\mathfrak l$, $
H^*(\phi_f[-1]\psi_{\mathfrak l}[-1]\Fdot)_\0$ is isomorphic to $H^*(\psi_{\mathfrak l}[-1]\phi_f[-1]\Fdot)_\0$.
\end{enumerate}
\end{cor}

\smallskip

\begin{rem}  Suppose that $X=\U$ and $\Fdot=\Z_\U^\bullet$ throughout this remark.

Then, $\hyp^{*}(F_{\mathfrak l, \mathbf 0} ; \Fdot)$ has the cohomology of a point and $\hyp^{*}(F_{{\mathfrak l}_{|_{V(f)}}, \mathbf 0} ; \Fdot)$ is the cohomology of the complex link of $V(f)$ at $\0$. Hence, the equivalence of Item 2 (or Item 1) and Item 5 of \corref{cor:emptyequiv2} is a generalization of the well-known result that the cohomology of the complex link of $V(f)$ at $\0$ is isomorphic to that of a point if and only if the relative polar curve is empty (in a neighborhood of $\0$).

Suppose now that we also have that $\dim_\0\Sigma f=1$. Then, for generic $\mathfrak l$, $f_{|_{V(\mathfrak l)}}$ has an isolated critical point at the origin, and $H^*(\phi_f[-1]\psi_{\mathfrak l}[-1]\Fdot)_\0$ is isomorphic to the (shifted) reduced integral cohomology of the Milnor fiber at the origin of $f_{|_{V(\mathfrak l)}}$. On the other hand, $H^*(\psi_{\mathfrak l}[-1]\phi_f[-1]\Fdot)_\0$ is isomorphic to the direct sum of the  reduced integral cohomologies of the Milnor fibers of $f_{|_{V(\mathfrak l-t)}}$, where  the sum is over all points $p\in\stackrel{\circ}{B}_\epsilon\cap\Sigma f\cap V(\mathfrak l - t)$ for $0<|t|\ll\epsilon\ll 1$. Therefore, the equivalence of Item 2 (or Item 1) and Item 6 of \corref{cor:emptyequiv2} is a generalization of the well-known result that the Milnor number of a generic hyperplane slice equals the sum of the Milnor numbers in a nearby hyperplane slice if and only if the relative polar curve is empty (in a neighborhood of $\0$).
\end{rem}

\section{An Application to Thom's $a_f$ Condition}\label{sec:afappl}

By combining our results from \cite{vanaf} with \thmref{thm:main1}, we can relate the polar curve to Thom's $a_f$ condition. Essentially what we prove below, in \thmref{thm:main3}, is that, if a stratification satisfies the $a_f$ condition, except perhaps at a point $p$ on a $1$-dimensional stratum, then the stratification satisfies the $a_f$ condition at $p$ if and only if, for some affine linear form $\mathfrak l$,  the polar curve of $(f, \mathfrak l)$ at $p$ is empty.

However, we do not actually need to start with a stratification, for we do not need the condition of the frontier. Also, of course, we want such a result with respect to a complex of sheaves. So, we need to make a number of preliminary definitions before we can state and prove our precise result.

\smallskip

Suppose that $M$ and $N$ are complex submanifolds of $\U$.

\begin{defn} The {\bf pair $(M, N)$ satisfies Thom's $a_{\tilde f}$ condition at a point $x\in N$} if and only if there is an inclusion, of fibers over $x$, $\big(\overline{T^*_{\tilde f_{|_M}}\U}\big)_x\subseteq \big(T^*_{\tilde f_{|_N}}\U\big)_x$.

The {\bf pair $(M, N)$ satisfies Thom's $a_{\tilde f}$ condition} if and only if it satisfies the $a_{\tilde f}$ condition at each point $x\in N$.
\end{defn}

\begin{rem} Note that if $\tilde f$ is a locally constant function, then the $a_{\tilde f}$ condition reduces to condition (a) of Whitney.
\end{rem}

The $a_{\tilde f}$ is condition is important for several reasons. First, it is an hypothesis of Thom's second isotopy lemma; see \cite{mather}. Second, the $a_{\tilde f}$ condition, and the existence of stratifications in which all pairs of strata satisfy the $a_{\tilde f}$ condition, is essential in arguments such as that used by L\^e in \cite{relmono} to prove that Milnor fibrations exist even when the domain is an arbitrarily singular space. Third, the $a_{\tilde f}$ condition is closely related to constancy of the Milnor number in families of isolated hypersurface singularities; see \cite{lesaito}.

There are at least two important general results about the $a_{\tilde f}$ condition: the above-mentioned existence of $a_{\tilde f}$ stratifications, proved first by Hironaka in \cite{hironakastratflat} and then in a different manner by Hamm and L\^e, following an argument of F. Pham,  in Theorem 1.2.1 of  \cite{hammlezariski}, and the fact that Whitney stratifications in which $V(\tilde f):=\tilde f^{-1}(0)$ is a union of strata  are $a_{\tilde f}$ stratifications, proved independently by Parusi\'nski  in \cite{parusw_f}, and Brian\c con, P. Maisonobe, and M. Merle in \cite{bmm}.

\smallskip

We can easily prove the following:

\begin{prop}\label{prop:easyaf} Suppose that $p\not\in\overline{\Sigma_{\Fdot}f}$. Let $M$ and $N$ be analytic submanifolds of $\U$ such that $\overline{T^*_M\U}$ is an irreducible component of $SS(\Fdot)$, and such that $(M, N)$ satisfies Whitney's condition (a) at a point $p\in N$. 

Then, $\big(\overline{T^*_{f_{|_M}}\U}\big)_p\subseteq (T^*_N\U)_p+<d_p\tilde f>$. In particular, if $\tilde f$ is locally constant on $N$ at $p$ or if $p\not\in\Sigma(\tilde f_{|_{N}})$, then  $(M, N)$ satisfies Thom's $a_{\tilde f}$ condition at $p$.
\end{prop}
\begin{proof} By replacing $\U$ with a small neighborhood of $p$, we may assume that $N$ is connected and closed in $\U$. Note that, as $\overline{T^*_M\U}$ is an irreducible component of $SS(\Fdot)$, $M$ is connected, $M\subseteq X$, and $\tilde f_{|_M}=f_{|_M}$.

There are 2 cases to consider.

\smallskip

\noindent Case 1: $\tilde f$ is constant on $M$.

By Whitney's condition (a),
$$
\big(\overline{T^*_{f_{|_M}}\U}\big)_p=\big(\overline{T^*_{M}\U}\big)_p\subseteq (T^*_N\U)_p\subseteq \big(T^*_N\U\big)_p+<d_p\tilde f>.
$$

\smallskip

\noindent Case 2: $\tilde f$ is not constant on $M$.

Since $p\not\in\overline{\Sigma_{\Fdot}f}$ and $\overline{T^*_M\U}$ is an irreducible component of $SS(\Fdot)$, \lemref{lem:lem} tells us that $\big(\overline{T^*_{f_{|_M}}\U}\big)_p = \big(\overline{T^*_{M}\U}\big)_p +<d_p\tilde f>$, and so, using Whitney's condition (a) again, we have
$$
\big(\overline{T^*_{f_{|_M}}\U}\big)_p = \big(\overline{T^*_{M}\U}\big)_p +<d_p\tilde f> \subseteq (T^*_N\U)_p+<d_p\tilde f>.
$$
\end{proof}

\begin{defn}\label{def:part} A collection $\W$ of subsets of $X$ is a {\bf (complex analytic) partition} of $X$ if and only if  $\W$ is a locally finite disjoint collection of analytic submanifolds of $\U$, which we call {\it strata}, whose union is all of $X$, and such that, for each stratum $W\in\W$, $\overline{W}$ and $\overline{W}-W$ are closed complex analytic subsets of $X$.

\smallskip

In this paper, we assume that all of the strata of a partition are connected.

\smallskip

A partition $\W$ is a {\bf stratification} if and only if it satisfies the condition of the frontier, i.e., for all $W\in\W$, $\overline W$ is a union of elements of $\W$.
\end{defn}

Below, we extend our earlier definition of $\Fdot$-visible strata to the case of a partition which may not satisfy Whitney conditions.

\begin{defn}\label{def:fdotpart}
A partition $\W$ of $X$ is an {\bf $\Fdot$-partition} provided that 
$$
SS(\Fdot)\subseteq \bigcup_{W\in\W}\overline{T^*_{{}_{W}}\U}.
$$

If $\W$ is an {\bf $\Fdot$-partition}, then a stratum $W\in\W$ is {\bf $\Fdot$-visible} if and only if $\overline{T^*_{{}_{W}}\U}\subseteq SS(\Fdot)$. We let $\W(\Fdot):= \{W\in\W\ |\ W\ {\rm is}\ \Fdot{\text -visible}\}$.
\end{defn}

\begin{rem} The reader should understand that the point of an $\Fdot$-partition $\W$ is that, for each $\Fdot$-visible stratum $S$ in $\strat$, there exists a unique $W\in\W$ such that $\overline{S}=\overline{W}$ and, hence, $\overline{T^*_{{}_{S}}\U}=\overline{T^*_{{}_{W}}\U}$. It follows at once from this, and the definition of $\Fdot$-visible strata of $\W$, that, if $\W$ is an $\Fdot$-partition, then
$$SS(\Fdot)\ =\ \bigcup_{W\in \W(\Fdot)}\overline{T^*_{{}_{W}}\U}.$$
\end{rem}

We can now give our result on Thom's $a_f$ condition and the relative polar curve.

\begin{thm}\label{thm:main3} Suppose that
\begin{enumerate}
\item[a.] $\W$ is an $\Fdot$-partition of $X$;
\item[b.] $\W^\prime$ is a Whitney (a) partition of $V(f)$;
\item[c.] $\0\in T\in\W^\prime$ and $\dim T=1$;
\item[d.] for all $W\in\W(\Fdot)$ such that $W\not\subseteq V(f)$, for all $W^\prime\in\W^\prime$ such that $W^\prime\neq T$, $(W, W^\prime)$ satisfies the $a_f$ condition.
\end{enumerate}

\noindent Then, the condition:

\noindent $(\dagger)$ for all $W\in\W(\Fdot)$ such that $W\not\subseteq V(f)$, $(W, T)$ satisfies the $a_f$ condition 

\noindent is equivalent to all of the conditions in \corref{cor:emptyequiv2}; in particular, it is equivalent to: there exists a non-zero linear form $\mathfrak l$ such that $\0\not\in\big|\Gamma_{f, \mathfrak l}(\Fdot)\big|$.
\end{thm}
\begin{proof} Let $\W^\prime_T:=\{W^\prime\in\W^\prime\ |\ W^\prime\neq T\}$. By Corollary 3.9 of \cite{vanaf}, our hypotheses imply that, if $p\in\W^\prime\in\W^\prime_T$, and ${\mathfrak l}^\prime$ is a non-zero linear form such that $V({\mathfrak l}^\prime-{\mathfrak l}^\prime(p))$ transversely intersects $W^\prime$ at $p$, then $p\not\in\supp\phi_{{\mathfrak l}^\prime-{\mathfrak l}^\prime(p)}[-1]\psi_f[-1]\Fdot$.

\medskip

Assume $(\dagger)$. By Corollary 3.9 of \cite{vanaf}, $\psi_f[-1]\Fdot$ is $\phi$-constructible with respect to $\W^\prime$. In particular, this implies that, if $\mathfrak l$ is a non-zero linear form such that $V(\mathfrak l)$ transversely intersects $T$ at $\0$, then $\0\not\in\supp\phi_{\mathfrak l}[-1]\psi_f[-1]\Fdot$. This implies Item 3 of  \corref{cor:emptyequiv2}.

Assume Item 3 of  \corref{cor:emptyequiv2}. Fix a non-zero $\mathfrak l$ such that $V(\mathfrak l)$ transversely intersects $T$ at $\0$ and $H^*(\phi_{\mathfrak l}[-1]\psi_f[-1]\Fdot)_\0=0$. As $\W^\prime$ satisfies Whitney (a), $V(\mathfrak l)$ transversely intersects all of the strata of $\W^\prime$ in a neighborhood of $\0$. Thus, by the first paragraph of the proof, $\0\not\in \supp\phi_{\mathfrak l}[-1]\psi_f[-1]\Fdot$. As $T$ is $1$-dimensional, this, together with the first paragraph of the proof, implies that $\psi_f[-1]\Fdot$ is weakly $\phi$-constructible with respect to $\W^\prime$. By Corollary 3.9 of \cite{vanaf}, this implies $(\dagger)$.
\end{proof}

\section{Families of Isolated Critical Points}\label{sec:families}

In \cite{critpts} and \cite{numcontrol}, we discussed continuous families of constructible complexes of sheaves. We wish to revisit our results in those works, and show how the proofs and results can be greatly improved by using the main theorems of this paper. 

We continue to let $\tilde f$ and $\tilde g$ be analytic functions from $\U$ to $\C$, and we let $f$ and $g$ be their respective restrictions to the analytic space $X$. We continue with $\Fdot$ being a bounded, constructible complex of sheaves of $R$-modules on $X$.

\smallskip

Throughout this section, we consider that $\Fdot$ and $g$ define a family (a {\it $g$-family}\/) of constructible complexes by setting $\Fdot_a:=\Fdot_{|_{V(g-a)}}[-1]$, for each $a\in\C$. We also consider the family of functions $f_a:=f_{|_{V(g-a)}}$. As in \cite{critpts} and \cite{numcontrol}, we make the following definition.

\begin{defn}\label{def:continuous} The family $\Fdot_a$ is {\bf continuous at $a=a_0$} if and only if $\phi_{g-a_0}[-1]\Fdot=0$. The family $\Fdot_a$ is {\bf continuous at $p\in X$} if and only if there exists an open neighborhood of $p$ in which $\phi_{g-g(p)}[-1]\Fdot=0$.
\end{defn}

\smallskip

The following is an additivity/upper-semicontinuity result. Recall our definition of the critical locus of a function relative to a complex of sheaves from \defref{def:critlocus}.

\smallskip

\begin{prop}\label{prop:additivity} Suppose that the family $\Fdot_a$ is continuous at $\0$, and that $\dim_0\overline{\Sigma_{\Fdot_0}f_0}\leq 0$. Then, $\dim_0\overline{\Sigma_{\Fdot}f}\leq 1$, and there exists an open neighborhood $\Omega$ of $\0$ in $\U$ and  $\delta>0$ such that, if $|a|<\delta$, then $\Omega\cap\overline{\Sigma_{\Fdot_a}f_a}$ is either empty or consists of a finite number of points and, for all $k$,
$$
H^k(\phi_{f_0}[-1]\Fdot_0)_\0\cong\bigoplus_{p\in\Omega\cap\overline{\Sigma_{\Fdot_a}f_a}}H^k(\phi_{f_a}[-1]\Fdot_a)_p.
$$
\end{prop}
\begin{proof} As $\phi_g[-1]\Fdot=0$ in a neighborhood of $\0$, $\phi_{f_0}[-1]\Fdot_0=\phi_f[-1]\psi_g[-1]\Fdot$ near $\0$, and so the assumption that $\dim_\0\overline{\Sigma_{\Fdot_0}f_0}\leq 0$ is equivalent to $\dim_\0\supp \phi_f[-1]\psi_g[-1]\Fdot\leq 0$. By \thmref{thm:main1}, this implies that $\dim_\0V(g)\cap|\Gamma_{g, \tilde f}(\Fdot)|\leq 0$ and 
$$
H^k(\phi_{f_0}[-1]\Fdot_0)_\0\cong H^k(\phi_f[-1]\psi_g[-1]\Fdot)_\0\cong  \big(\big(\Gamma^1_{g, \tilde f}(\Fdot)\big)^k\odot V(g)\big)_\0.
$$
Hence, there exists an open neighborhood $\Omega$ of $\0$ in $\U$ and  $\delta>0$ such that, if $|a|<\delta$, then $\Omega\cap V(g-a)\cap|\Gamma_{g, \tilde f}(\Fdot)|$ is either empty or consists of a finite number of points, and 
$$
H^k(\phi_{f_0}[-1]\Fdot_0)_\0\cong \bigoplus_{p\in\Omega\cap V(g-a)\cap|\Gamma_{g, \tilde f}(\Fdot)|}\big(\big(\Gamma^1_{g, \tilde f}(\Fdot)\big)^k\odot V(g-a)\big)_p\cong \bigoplus_{p\in\Omega\cap\overline{\Sigma_{\Fdot_a}f_a}}H^k(\phi_{f_a}[-1]\Fdot_a)_p,
$$
where the last equality follows again from \thmref{thm:main1}, and using that, near $\0$, $\overline{\Sigma_{\Fdot}g}\subseteq V(g)$.
\end{proof}

\bigskip

There is the trivial case:

\begin{prop}\label{prop:trivfam} Suppose that the family $\Fdot_a$ is continuous at $\0$, and that $\0\not\in\overline{\Sigma_{\Fdot_0}f_0}$. Then, 
\begin{enumerate}
\item $\0\not\in \overline{\Sigma_{\Fdot}f}$;
\item near $\0$, $\phi_g[-1]\big(\psi_f[-1]\Fdot\big)=\phi_g[-1]\big(\Fdot_{|_{V(f)}}[-1]\big)=0$, i.e.,  $\psi_f[-1]\Fdot$, and $\Fdot_{|_{V(f)}}[-1]$ are $g$-continuous families at $\0$; and
\item if $M$ and $N$ are analytic submanifolds of $\U$ such that $\overline{T^*_M\U}$ is an irreducible component of $SS(\Fdot)$, such that $(M, N)$ satisfies Whitney's condition (a) at $\0$, and such that  $N\subseteq V(f)$, or $\0\not\in\Sigma (f_{|_{N}})$, then $(M, N)$ satisfies Thom's $a_{\tilde f}$ condition at $\0$.
\end{enumerate}
\end{prop}
\begin{proof}

\noindent Proof of Item 1:

By \propref{prop:sigma},
$$\overline{\Sigma_{\Fdot}g}\ \cup\ \big|\Gamma_{g, \tilde f}(\Fdot)\big|={\overline{\Sigma}}_{\Fdot}(\tilde f, \tilde g)=\overline{\Sigma_{\Fdot}f}\ \cup\ \big|\Gamma_{f, \tilde g}(\Fdot)\big|.
$$
As $\0\not\in\overline{\Sigma_{\Fdot}g}$, there is an equality of sets $\big|\Gamma_{g, \tilde f}(\Fdot)\big| =\overline{\Sigma_{\Fdot}f}\ \cup\ \big|\Gamma_{f, \tilde g}(\Fdot)\big|$ near $\0$. As $\0\not\in\overline{\Sigma_{\Fdot_0}f_0}$, \propref{prop:additivity} implies that $\0\not\in\big|\Gamma_{g, \tilde f}(\Fdot)\big|$. Item 1 follows.

\medskip

\noindent Proof of Item 2:

It follows from the paragraph above that $\0\not\in \big|\Gamma_{f, \tilde g}(\Fdot)\big|$. As $\big|\Gamma_{f, \tilde g}(\Fdot)\big|$ is closed, we may apply \thmref{thm:main1} at each point $p$ near the origin in $V(f, g)$ to conclude that, near $\0$, $\phi_g[-1]\big(\psi_f[-1]\Fdot\big)=0$. By Item 1, $\phi_f[-1]\Fdot=0$ near the origin, and so $\psi_f[-1]\Fdot= \Fdot_{|_{V(f)}}[-1]$ near $\0$. Item 2 follows.

\medskip

\noindent Proof of Item 3:

This follows immediately from Item 1 and \propref{prop:easyaf}.
\end{proof}

\bigskip

The following theorem contains generalizations of well-known properties/results that hold in the classic case of families of isolated critical points of functions on affine space, including a generalization of the main result of L\^e and Saito from \cite{lesaito}. 

\begin{thm}\label{thm:fam} Suppose that the family $\Fdot_a$ is continuous at $\0$, and that $\dim_0\overline{\Sigma_{\Fdot_0}f_0}\leq 0$. 

Let $C\subseteq V(f)$ be a locally irreducible curve which contains $\0$ such that, for all $p\in C$ near $\0$, the isomorphism-type of $H^*(\phi_{f_{g(p)}}[-1]\Fdot_{g(p)})_p$ is non-zero and equal to $H^*(\phi_{f_{0}}[-1]\Fdot_{0})_\0$. Then, in a neighborhood of the origin,

\begin{enumerate}
\item $C$ is smooth;

\item $V(\tilde g)$ is smooth, and transversely intersects $C$;

\item $\0\not\in\big|\Gamma_{f, \tilde g}(\Fdot)\big|$;

\item $C=\overline{\Sigma_{\Fdot}f} = \big|\Gamma_{g, \tilde f}(\Fdot)\big|$;

\item  $\phi_g[-1]\big(\phi_f[-1]\Fdot)=0$, $\phi_g[-1]\big(\psi_f[-1]\Fdot\big)=0$, and $\phi_g[-1]\big(\Fdot_{|_{V(f)}}[-1]\big)=0$, i.e.,  $\phi_f[-1]\Fdot$, $\psi_f[-1]\Fdot$, and $\Fdot_{|_{V(f)}}[-1]$ are $g$-continuous families at $\0$;

\item $SS(\phi_f[-1]\Fdot)=T^*_C\U$;

\item  for all $k$, the sheaf cohomology $\mathbf H^{k-1}(\phi_f[-1]\Fdot)$ is constant on $C$ and, for all $p\in C$, there is an isomorphism of stalk cohomology $H^{k-1}(\phi_f[-1]\Fdot)_p\cong H^{k}(\phi_{f_{0}}[-1]\Fdot_{0})_\0$ ;

\item for all $p\in C$, the Milnor monodromy automorphisms 
$$\widetilde T^*_{f_{g(p)}, p}:H^*(\phi_{f_{g(p)}}[-1]\Fdot_{g(p)})_p\rightarrow H^*(\phi_{f_{g(p)}}[-1]\Fdot_{g(p)})_p$$
 are isomorphic to $\widetilde T^*_{f_0, \0}$, i.e., there are isomorphisms between $H^*(\phi_{f_{g(p)}}[-1]\Fdot_{g(p)})_p$ and $H^*(\phi_{f_{0}}[-1]\Fdot_{0})_\0$ which commute with the monodromy automorphisms; and

\item  if $M$ is an analytic submanifold of $\U$ such that $\overline{T^*_M\U}$ is an irreducible component of $SS(\Fdot)$ and $(M, C)$ satisfies Whitney's condition $(a)$ at $\0$, then $(M, C)$ satisfies Thom's $a_f$ condition at $\0$.
\end{enumerate}
\end{thm}
\begin{proof} As we saw in \propref{prop:additivity}, the assumptions imply that there exists an open neighborhood $\Omega$ of $\0$ in $\U$ and  $\delta>0$ such that, if $|a|<\delta$,
$$
H^k(\phi_{f_0}[-1]\Fdot_0)_\0\cong \bigoplus_{p\in\Omega\cap V(g-a)\cap|\Gamma_{g, \tilde f}(\Fdot)|}\big(\big(\Gamma^1_{g, \tilde f}(\Fdot)\big)^k\odot V(g-a)\big)_p\cong \bigoplus_{p\in\Omega\cap\overline{\Sigma_{\Fdot_a}f_a}}H^k(\phi_{f_a}[-1]\Fdot_a)_p.
$$
Hence, the constancy of the isomorphism-type of $H^*(\phi_{f_{g(p)}}[-1]\Fdot_{g(p)})_p$ along $C$ implies that, near $\0$,  $C=\big|\Gamma_{g, \tilde f}(\Fdot)\big|$, $C$ is smooth at $\0$, and $V(\tilde g)$ is smooth at $\0$, and transversely intersects $C$ at $\0$. This proves Items 1 and 2.

As we saw in the proof of \propref{prop:trivfam}, near $\0$, $\big|\Gamma_{g, \tilde f}(\Fdot)\big| =\overline{\Sigma_{\Fdot}f}\ \cup\ \big|\Gamma_{f, \tilde g}(\Fdot)\big|$; as we saw above, this equals the irreducible curve $C$. As $\dim_\0 V(g)\cap \big|\Gamma_{g, \tilde f}(\Fdot)\big|= 0$, $\dim_\0 V(g)\cap \big|\Gamma_{f, \tilde g}(\Fdot)\big|\leq 0$. By \lemref{lem:phipsi}, this implies that $\dim_\0 V(f)\cap \big|\Gamma_{f, \tilde g}(\Fdot)\big|\leq 0$. As every component of $\big|\Gamma_{f, \tilde g}(\Fdot)\big|$ has dimension at least one, and as $\big|\Gamma_{f, \tilde g}(\Fdot)\big|$ is contained in $C\subseteq V(f)$ near $\0$, we conclude that $\0\not\in\big|\Gamma_{f, \tilde g}(\Fdot)\big|$; this is Item 3. As $\big|\Gamma_{g, \tilde f}(\Fdot)\big| =\overline{\Sigma_{\Fdot}f}\ \cup\ \big|\Gamma_{f, \tilde g}(\Fdot)\big|=C$ near $\0$, we conclude Item 4.

From Item 3 of this theorem and \thmref{thm:main1}, we conclude that the stalk cohomology at the origin of $\phi_g[-1]\big(\psi_f[-1]\Fdot\big)$ is equal to zero. As we may apply this argument at each point near $\0$, we conclude that $\phi_g[-1]\big(\psi_f[-1]\Fdot\big)=0$ near $\0$. From Item 3 of this theorem and Item 2 of \corref{cor:main2}, we find that $\hyp^*(F_{g, \0}, F_{g_{|_{V(f)}}, \0}; \Fdot)=0$. As $\hyp^*(B_\epsilon, F_{g, \0}; \Fdot)=0$, we conclude that $\hyp^*(B_\epsilon\cap V(f), F_{g_{|_{V(f)}}, \0}; \Fdot)=0$, i.e., that the stalk cohomology at $\0$ of $\phi_g[-1]\big(\Fdot_{|_{V(f)}}[-1]\big)$ is zero. Since we may apply this argument at each point near $\0$, we conclude that $\phi_g[-1]\big(\Fdot_{|_{V(f)}}[-1]\big)=0$ near the origin. As $\phi_g[-1]\big(\psi_f[-1]\Fdot\big)=0$ and $\phi_g[-1]\big(\Fdot_{|_{V(f)}}[-1]\big)=0$ near $\0$, the distinguished triangle relating the nearby and vanishing cycles implies that  $\phi_g[-1]\big(\phi_f[-1]\Fdot\big)=0$ near $\0$. This proves Item 5.

\smallskip

Item 6 follows at once from Items 4 and 5, together with the fact that $d_\0\tilde g\neq 0$ by Item 2. Item 7 follows immediately from Item 4 and the fact, from Item 5, that $\phi_g[-1]\big(\phi_f[-1]\Fdot\big)=0$ near $\0$. 

\smallskip

Items 3 and 5 tell us that, in  \thmref{thm:main2}, $\Adot:=R(T^\epsilon_{\delta, \rho})_*(\Fdot)^\epsilon_{\delta, \rho}$ is complex analytically constructible with respect to the stratification given by $\{{\stackrel{\circ}{\D}}_{{}_\delta}\times {\stackrel{\circ}{\D}}_{{}_\rho}-{\stackrel{\circ}{\D}}_{{}_\delta}\times\{0\},\ {\stackrel{\circ}{\D}}_{{}_\delta}\times\{0\}\}$. Using $(u, v)$ for coordinates on ${\stackrel{\circ}{\D}}_{{}_\delta}\times {\stackrel{\circ}{\D}}_{{}_\rho}$, it is trivial that the Milnor monodromy on $H^*(\phi_v[-1](\Adot_{|_{V(u-a)}})_{(a, 0)}$ is constant (up to isomorphism) for small $|a|$. As $C=\overline{\Sigma_{\Fdot}f}=\supp\phi_f[-1]\Fdot$ near $\0$, Item 8 follows.

\smallskip

Let $M$ be as in Item 9, and let $Y:=\overline{M}$. If $\0\not\in Y$ or if $Y=C$ near $\0$, then Item 9 follows trivially. So, suppose that $\0\in Y$ and $Y\neq C$ near $\0$. As $\overline{T^*_M\U}$ is an irreducible component of $SS(\Fdot)$, $Y$ is irreducible. If $f$ were constant on $Y$, then, as $\0\in Y$, $Y$ would have to be contained in $V(f)$, and so we would have $Y\subseteq\overline{\Sigma_{\Fdot}f}$; this is impossible by Item 4. 

Now, let $E$ denote the exceptional divisor in the blow-up ${\operatorname{Bl}}_{\operatorname{im}d\tilde f}(\overline{T^*_M\U})\subseteq T^*\U\times\Proj^n$. We identify $T^*\U\times\Proj^n$ with $\U\times\C^{n+1}\times\Proj^n$, and let $\sigma:\U\times\C^{n+1}\times\Proj^n\rightarrow \U\times\Proj^n$ denote the projection. Over a neighborhood of the origin, Item 6 tells us that $SS(\phi_f[-1]\Fdot)=T^*_C\U$. Now, by Theorem 3.4 of \cite{singenrich}, $\sigma(E)\cap(\{\0\}\times\Proj^n)\subseteq \Proj(T^*_C\U)_\0$. By Proposition 4.3 of \cite{pervcohovan}, this implies that $(M, C)$ satisfies Thom's $a_f$ condition at $\0$.
\end{proof}

\begin{rem} The reader should compare the statements and the proofs from \thmref{thm:fam} with our related results in \cite{critpts} and \cite{numcontrol}. Not only are our current results more general and stronger, the proofs are vastly easier.

We also remark that Item 5 of \thmref{thm:fam} is useful for inductions. For instance, in proving results for families of local complete intersections with isolated singularities.
\end{rem}

\bigskip

Recall the definition of the graded, enriched characteristic cycle from \defref{def:gecc}. Let $SS^k(\Fdot)$ equal the underlying set $|\gecc^k(\Fdot)|$. This means that $SS^k(\Fdot)=\bigcup \overline{T^*_S\U}$, where the union is over those strata $S\in\strat$ such that $H^{k-d_S}(\mathbb N_S,\mathbb L_S; \Fdot)\neq 0$.
We also set $\Sigma^k_{\Fdot} f:=\{p\in X\ |\ H^k(\phi_{f-f(p)}[-1]\Fdot)_p\neq 0\}$.

In Proposition 2.5 of \cite{singenrich}, we explained how the perverse cohomology operator ${}^{\mu}\hskip -.02in H^k$ works with graded, enriched characteristic cycles; we showed there is an equality of enriched cycles given by 
$$
\gecc^0\big({}^{\mu}\hskip -.02in H^k(\Fdot)\big) = \gecc^k(\Fdot).
$$
This is a consequence of the fact that the perverse cohomology operator commutes with the $[-1]$-shifted vanishing cycle operator, and that, if $p$ is an isolated point in the support of a complex $\Adot$, then $H^j\big({}^{\mu}\hskip -.02in H^k(\Adot)\big)_p$ is isomorphic to $H^k(\Adot)_p$ when $j=0$, and is zero if $j\neq 0$.

The following corollary, which tells us that \thmref{thm:fam} applies one degree at a time, follows immediately. Note, however, that the conditions in the first line are {\bf not} conditions on a fixed degree.

\begin{cor}\label{cor:fam} Suppose that the family $\Fdot_a$ is continuous at $\0$, and that $\dim_0\overline{\Sigma_{\Fdot_0}f_0}\leq 0$. Fix $k\in \Z$.

Let $C\subseteq V(f)$ be a locally irreducible curve which contains $\0$ such that, for all $p\in C$ near $\0$, the isomorphism-type of $H^k(\phi_{f_{g(p)}}[-1]\Fdot_{g(p)})_p$ is non-zero and equal to $H^k(\phi_{f_{0}}[-1]\Fdot_{0})_\0$. Then, in a neighborhood of the origin,

\begin{enumerate}
\item $C$ is smooth;

\item $V(\tilde g)$ is smooth, and transversely intersects $C$;

\item $\0\not\in\big|\big(\Gamma^1_{f, \tilde g}(\Fdot)\big)^k\big|$;

\item $C=\overline{\Sigma^{k-1}_{\Fdot}f} = \big|\big(\Gamma^1_{g, \tilde f}(\Fdot)\big)^k\big|$;

\item $SS^{k}(\phi_f[-1]\Fdot)=T^*_C\U$;

\item the sheaf cohomology $\mathbf H^{k-1}(\phi_f[-1]\Fdot)$ is constant on $C$ and, for all $p\in C$, there is an isomorphism of stalk cohomology $H^{k-1}(\phi_f[-1]\Fdot)_p\cong H^{k}(\phi_{f_{0}}[-1]\Fdot_{0})_\0$ ;

\item for all $p\in C$, the Milnor monodromy automorphisms 
$$\widetilde T^k_{f_{g(p)}, p}:H^k(\phi_{f_{g(p)}}[-1]\Fdot_{g(p)})_p\rightarrow H^k(\phi_{f_{g(p)}}[-1]\Fdot_{g(p)})_p$$
 are isomorphic to $\widetilde T^k_{f_0, \0}$, i.e., there are isomorphisms between $H^k(\phi_{f_{g(p)}}[-1]\Fdot_{g(p)})_p$ and $H^k(\phi_{f_{0}}[-1]\Fdot_{0})_\0$ which commute with the monodromy automorphisms; and

\item  if $M$ is an analytic submanifold of $\U$ such that $\overline{T^*_M\U}$ is an irreducible component of $SS^k(\Fdot)$ and $(M, C)$ satisfies Whitney's condition $(a)$ at $\0$, then $(M, C)$ satisfies Thom's $a_f$ condition at $\0$.
\end{enumerate}
\end{cor}

\bibliographystyle{plain}
\bibliography{Masseybib}
%\printindex
\end{document}